\documentclass[oneside,reqno]{amsart}
\usepackage{amssymb,latexsym} 
\usepackage{verbatim} 
\usepackage{eucal} 
\usepackage[normalem]{ulem} 

\theoremstyle{plain}
\newtheorem{theorem}{Theorem}[subsection]
\newtheorem{lemma}[theorem]{Lemma}
\newtheorem{corollary}[theorem]{Corollary}
\newtheorem{proposition}[theorem]{Proposition}

\theoremstyle{definition}

\newtheorem{definition}[theorem]{Definition}
\newtheorem{notation}[theorem]{Notation}

\newtheorem{example}[theorem]{Example}
\newtheorem{remark}[theorem]{Remark}

\newcommand\catA{{\mathfrak A_\text{cat}}}
\newcommand\catB{{\mathfrak B_\text{cat}}}
\newcommand\classA{{\mathfrak A(\bar\Lambda)}}

\newcommand\classB{{\mathfrak B(\Lambda,\Gamma)}}

\newcommand\classBfi{\mathfrak B(\Lm,\textrm{f\hskip1.5pt i})}

\newcommand\cA{\mathcal A}
\newcommand\cB{\mathcal B}

\newcommand\cI{{\mathcal I}}

\newcommand\cL{{\mathcal L}}
\newcommand\cM{{\mathcal M}}

\newcommand\cS{{\mathcal S}}

\newcommand{\blm}{{\bar\lambda}}
\newcommand{\bLm}{{\bar\Lambda}}
\newcommand{\bmu}{{\bar\mu}}
\newcommand{\gm}{\gamma}
\newcommand{\Gm}{\Gamma}
\newcommand{\lm}{\lambda}
\newcommand{\Lm}{\Lambda}
\newcommand{\ph}{\varphi}
\newcommand{\sg}{\sigma}

\newcommand\bbZ{\mathbb Z}
\newcommand\Zn{{\mathbb Z^n}}
\newcommand\Ztwo{{\mathbb Z^2}}

\newcommand\boldkb{\bar{\mathbf{k}}}
\newcommand\boldl{{\boldsymbol{\ell}}}
\newcommand\boldlb{\,\,{\bar{\!\!\boldsymbol\ell}}}

\newcommand\boldm{{\mathbf m}}
\newcommand\boldmp{{\mathbf m'}}

\newcommand\boldq{\mathbf q}

\newcommand\Achi{\cA_\chi}
\newcommand\algC{\Alg(C(\cB),k)}
\newcommand\andd{\quad \text{and} \quad}
\newcommand\blmsum{\bigoplus_{\bar \lambda\in\bar \Lambda}}
\newcommand\charG{\Hom(\Gamma,k^\times)}
\newcommand\dotchi{\cdot_\chi}

\newcommand\kerpi{\ker(\pi)}
\newcommand\kerrhoB{\ker(\rho)\cB}
\newcommand\ig{\simeq_{\text{ig}}}
\newcommand\isomLm{\simeq_\Lambda}
\newcommand\isombLm{\simeq_{\bar \Lambda}}
\newcommand\lmsum{\bigoplus_{\lambda\in\Lambda}}
\newcommand\loopA{L_\pi(\mathcal A)}
\newcommand\loopAp{L_\pi(\mathcal A')}
\newcommand\loopASig{L_\pi(\mathcal A,\Sigma)}

\newcommand\mloopA{M_{\boldm}(\cA,\sg_1,\dots,\sg_n)}
\newcommand\mloopAp{M_{\boldmp}(\cA',\sg'_1,\dots,\sg'_n)}

\newcommand\modell{{\bar{\ell}}}
\newcommand\order[1]{\left\vert #1 \right\vert}
\newcommand\ot{\otimes}
\newcommand\otk{\otimes_k}
\newcommand\projrho{p_\rho}
\newcommand \set[1]{\{ \, #1 \, \}}

\newcommand\unim{\cB/\ker(\rho)\cB}

\newcommand\unten{\cB \otimes_{C(\cB)} k_\rho}
\newcommand\suchthat{\mid}

\newcommand\Alg{\operatorname{Alg}}
\newcommand\Aut{\operatorname{Aut}}
\newcommand\charr{\operatorname{char}}
\newcommand\diag{\operatorname{diag}}
\newcommand\GL{\operatorname{GL}}

\newcommand\End{\operatorname{End}}
\newcommand\Hom{\operatorname{Hom}}
\newcommand\Mult{\operatorname{Mult}}
\newcommand\spann{\operatorname{span}}

\newcommand\supp{\operatorname{supp}}

\begin{document}

\title[Realization of graded-simple algebras]{Realization of graded-simple algebras as loop algebras}
\author{Bruce Allison}
\address[Bruce Allison]
{Department of Mathematical and Statistical Sciences\\ University of
Alberta\\Edmonton, Alberta, Canada T6G 2G1}
\email{ballison@ualberta.ca}
\author{Stephen Berman}
\address[Stephen Berman]
{Saskatoon, Saskatchewan, Canada}
\email{sberman@shaw.ca}
\author{John Faulkner}
\address[John Faulkner]
{Department of Mathematics\\
University of Virginia\\
Kerchof Hall, P.O.~Box 400137\\
Charlottesville VA 22904-4137 USA}
\email{jrf@virginia.edu}
\author{Arturo Pianzola}
\address[Arturo Pianzola]
{Department of Mathematical and Statistical Sciences\\ University of
Alberta\\Edmonton, Alberta, Canada T6G 2G1}
\email{a.pianzola@ualberta.ca} \dedicatory{Dedicated to the memory
of Gordon E.~Keller, 1939--2003}
\thanks{The authors Allison, Berman and Pianzola
gratefully acknowledge the support
of the Natural Sciences and Engineering Research
Council of Canada.  They also thank the University
of Virginia for its hospitality
during part of the work on this paper.}

\begin{abstract}
Multiloop algebras determined by $n$ commuting algebra automorphisms of finite order
are natural generalizations of the classical loop algebras
that are used to realize affine Kac-Moody Lie algebras.  In this
paper, we obtain necessary and sufficient conditions
for a $\bbZ^n$-graded algebra to be realized as a
multiloop algebra based on a finite dimensional
simple algebra over an algebraically closed field of characteristic 0.
We also obtain necessary and sufficient
conditions for two such multiloop algebras to be graded-isomorphic,
up to automorphism of the grading group.

We prove these facts as consequences
of corresponding results for a generalization of the multiloop construction.
This more general setting allows us to work
naturally and conveniently with arbitrary
grading groups and arbitrary base fields.
\end{abstract}
\date{\today}
\maketitle

\section{Introduction}
\label{sec:introduction}

This  paper studies the realization, or construction, of
graded-simple algebras as loop algebras.  Our  results are
quite general  and apply to algebras of any kind including both
Lie algebras and associative algebras. However, it  was a
very specific problem in the theory of infinite dimensional Lie
algebras that motivated our  work.

In V.~Kac's early work on infinite dimensional Lie algebras, he
showed that any affine Kac-Moody Lie algebra (or more precisely the
derived algebra modulo its centre of any affine Kac-Moody Lie
algebra) can be realized as the loop algebra of a finite order
automorphism of a finite dimensional  simple Lie algebra.
This fact  is of great importance in the theory of affine algebras.
Now extended affine Lie algebras (EALA's) are higher nullity
generalizations of affine Kac-Moody Lie algebras \cite{AABGP},
and so it
is natural to ask if any EALA (or more precisely the centreless core
of any EALA) can be realized as the multiloop algebra of a sequence
of commuting finite order automorphisms of a finite dimensional
 simple Lie algebra. Our  research on multiloop
realizations began with  this question.

The centreless core $\cL$ of an EALA, now also called a centreless
Lie torus  (\cite{N,Y2}),
is in particular graded by a finitely generated abelian group $\Lm$
of finite rank $n$, where $n$ is the nullity of the EALA.  Moreover,
as a $\Lm$-graded algebra, $\cL$ is graded-central-simple.
In fact, as our work progressed, it became clear that
the methods we were employing  applied  not only to centreless Lie tori
but also to graded-central-simple
algebras of any kind.

Indeed suppose that $\Lm$ is a free abelian group of finite rank
and $k$ is an algebraically closed field of characteristic 0.
We are able to show in Corollary \ref{cor:multreal}
that  a $\Lm$-graded algebra $\cB$ has a
multiloop realization based on a finite dimensional  simple
algebra if and only if $\cB$ is graded-central-simple,
$\cB$ is a finitely generated module over its centroid $C(\cB)$,
and the support $\Gm(\cB)$ of $C(\cB)$ has finite index in~$\Lm$.
This last condition on $\Gm(\cB)$ is redundant
in most cases of interest including Lie tori.
Thus, returning to our original problem, it follows that
a centreless Lie torus has a multiloop realization based on a finite dimensional
 simple algebra if and only if it is finitely generated as
a module over its centroid.

We also obtain a graded-isomorphism theorem for
two multiloop algebras based on a finite dimensional
simple algebra.  Roughly speaking
the theorem states that two such multiloop algebras are graded-isomorphic,
up to automorphism of the grading group,
if and only if the sequences of commuting automorphisms
that determine the algebras are in the same orbit under
 the natural action of $\GL_n(\bbZ)$.
(See Theorem \ref{thm:realmult}(ii) for a more precise statement.)

For most of the paper, we work in a very general setting.  We assume
that $\Lm$ is an arbitrary abelian group, $k$ is an arbitrary field,
 and study arbitrary graded-central-simple $\Lm$-graded algebras over $k$.
Rather than work with the multiloop algebra construction, which requires
primitive roots of unity in $k$, we work instead with a more general loop construction.
Given a fixed group epimorphism $\pi : \Lm \to \bLm$ with kernel $\Gm$, this general
loop construction
$L_\pi$ produces a $\Lm$-graded algebra $\loopA$ from
a $\bLm$-graded algebra $\cA$.  This general point of view
provides useful additional flexibility and simplicity even though
our main interest is in multiloop algebras.

To study the construction $L_\pi$,
it is convenient to introduce two classes of graded algebras:
$\classA$ is the class of all $\bLm$-graded algebras
that are central-simple as algebras, whereas
$\classB$ is the class of all $\Lm$-graded algebras
$\cB$ such that $\cB$ is graded-central-simple, $\Gm(\cB) = \Gm$,
and the centroid
$C(\cB)$ is split (isomorphic to the group algebra $k[\Gm]$).
If $\cA$ is in $\classA$, then $\loopA$ is in $\classB$.
Moreover the main result of the paper, which we call the Correspondence
Theorem (Theorem \ref{thm:correspondence}), states that
$L_\pi$ establishes a 1-1 correspondence  between similarity classes
of $\bLm$-graded algebras in $\classA$ and isomorphism classes
of $\Lm$-graded algebras in $\classB$.  (Similarity is discussed
in \S 6.3.)  To obtain the inverse of this correspondence
we construct from any $\cB\in \classB$  the quotient
algebra $\cB/\ker(\rho)\cB$ in $\classA$, where
$\rho: C(\cB) \to k$ is an arbitrary  unital algebra homomorphism.
We call this quotient algebra a central image of $\cB$.
The freedom to choose $\rho$ explains  why we work
with similarity classes in $\classA$ rather than
with graded-isomorphism classes.

All of the results mentioned above about multiloop algebras
are obtained as consequences of the Correspondence Theorem.

To conclude this introduction, we briefly describe the contents of the paper.
After a short preliminary section on graded algebras, we describe
in Section \ref{sec:loop} the general loop construction
and the multiloop construction.
In Section \ref{sec:graded-central-simple} we obtain some basic properties of the centroid
and graded-central-simple algebras, and in Section \ref{sec:gradedloop} we look at those
properties for loop algebras.  In Section \ref{sec:central-spec},
we study central images and similarity.

Section \ref{sec:correspondence} contains the Correspondence Theorem.
To illustrate that there is  interest in the case when $k$ is not necessarily
algebraically closed we include an example of the correspondence in the associative case.
In this example  $\classA$ contains many nonsplit finite dimensional
central-simple algebras, all of which   correspond to the same infinite dimensional
algebra in $\classB$, namely the quantum torus.

Section \ref{sec:mult} applies the Correspondence Theorem to obtain our
results about multiloop algebras.  Finally, in Section \ref{sec:tori}
we discuss applications of our results to three classes of algebras
that arise naturally as coordinate algebras in the study of EALA's.
These are the associative, alternative and Jordan tori.
We also briefly discuss the main application of this work to the study
of Lie tori, but a detailed discussion of this application will be
written in a sequel
to this paper.

The authors  thank  Karl-Hermann Neeb, Erhard Neher and Ottmar Loos
for helpful discussions regarding this work.

\section{Preliminaries}
\label{sec:prelim}

Throughout this work $k$ denotes an arbitrary field.
All algebras are assumed to be algebras (not necessarily associative or unital) over
$k$.  We also assume that $\Lm$ is an abelian group written additively.

In this section we recall some definitions and notation for $\Lm$-graded algebras.

\subsection{Definitions and notation}
\label{subsec:defn}

\begin{definition}
A  $\Lm$-\textit{graded algebra} is a pair $(\cB,\Sigma)$
consisting of an algebra $\cB$ together with a family $\Sigma =
\set{\cB^\lm}_{\lm\in\Lm}$ of subspaces of $\cB$ such that
$\cB = \lmsum \cB^\lm$ and $\cB^\lm \cB^\mu \subseteq \cB^{\lm+\mu}$
for $\lm,\mu\in\Lm$.  We call $\Sigma$ the $\Lm$-\textit{grading} of
$(\cB,\Sigma)$. We will usually suppress the symbol $\Sigma$ in the notation
and write the $\Lm$-graded algebra $(\cB,\Sigma)$ simply as $\cB$.
\end{definition}

\begin{example}
An important example of a $\Lm$-graded algebra is the group algebra
$k[\Lm] = \lmsum k z^\lm$
of $\Lm$, where the multiplication is given by $z^\lm z^\mu = z^{\lm + \mu}$
and the $\Lm$-grading is given by
$k[\Lm]^\lm = kz^\lm$ for $\lm\in \Lm$.
\end{example}

\begin{definition}
\label{def:isomorphism}
There are two notions of isomorphism that we will use for graded algebras.

(i) Suppose that $\cB$ and $\cB'$
are $\Lm$-graded algebras.  We say that $\cB$ and $\cB'$
are \textit{graded-isomorphic}, if there exists an algebra isomorphism
$\varphi: \cB \to \cB'$ such that $\varphi(\cB^\lm) = \cB'^\lm$ for $\lm\in \Lm$.
In that case we write $\cB \isomLm \cB'$.

(ii) Suppose that $\cB$ is a $\Lm$-graded algebra and $\cB'$ is a
$\Lm'$-graded algebra, where $\Lm'$ is another abelian group.  We
say that $\cB$ and $\cB'$ are \textit{graded-isomorphic up to
isomorphism of grading groups}, or more simply
\textit{isograded-isomorphic}, if there exists an algebra
isomorphism $\varphi: \cB \to \cB'$ and a group isomorphism $\nu :
\Lm \to \Lm'$ such that $\varphi(\cB^\lm) = \cB'^{\,\nu(\lm)}$ for
$\lm\in \Lm$. In that case we write $\cB \ig \cB'$.
\end{definition}

\begin{definition}
\label{def:regrade}
Let $\cB$ be a $\Lm$-graded algebra and let $\nu\in \Aut(\Lm)$.
We may \textit{regrade} $\cB$ using $\nu$ to obtain a new $\Lm$-graded
algebra $\cB_\nu$ as follows \cite[\S 1.1]{NVO1}:
As algebras $\cB_\nu$ and $\cB$ are the same, but the $\Lm$-grading
on $\cB_\nu$ is defined by $(\cB_\nu)^\lm = \cB^{\nu(\lm)}$.
\end{definition}

\begin{remark}
\label{rem:regrade}
Suppose that $\cB$ and $\cB'$ are $\Lm$-graded algebras.
Then $\cB \ig \cB'$ if and only if $\cB \isomLm \cB'_\nu$ for some $\nu\in \Aut(\Lm)$.
\end{remark}

\begin{notation}
If $\cB =\lmsum\cB^\lm$
is a $\Lm$-graded algebra over $k$, we use the notation
\[
\supp_\Lm(\cB) := \set{\lm\in\Lm \suchthat \cB^\lm \ne 0}
\]
for the \textit{$\Lm$-support} of  $\cB$.
We  denote the subgroup of $\Lm$ generated by $\supp_\Lm(\cB)$
as $\langle \supp_\Lm(\cB) \rangle$.
\end{notation}

\section{Loop algebras and  multiloop algebras}
\label{sec:loop}

In this section we introduce the main  object of our
study---the loop algebra $\loopA$.
The definition (and many of the results of this
paper) assumes only that $\Lm$ is an abelian group.

\subsection{The general definition}
\label{subsec:loop}
The  following definition is given in the associative case in
\cite[Proposition 1.2.2]{NVO2}.

\begin{definition}
\label{def:loop}
Suppose that $\pi : \Lm \to \bLm$ is a group epimorphism
of an abelian group $\Lm$ onto an abelian group $\bLm$.  We write
\[\blm = \pi(\lm)\]
for $\lm\in \Lm$.  Suppose that $\cA = \blmsum \cA^\blm$ is a  $\bLm$-graded algebra.
Then the tensor product $\cA\otk k[\Lm]$ is a $\Lm$-graded algebra over $k$, where
$(\cA\otk k[\Lm])^\lm = \cA\ot z^\lm$ for $\lm\in \Lm$.  We define
\[\loopA = \sum_{\lm\in\Lm} \cA^\blm \ot z^\lm\]
in $\cA\otk k[\Lm]$.  Then $\loopA$ is a $\Lm$-graded subalgebra of $\cA\otk k[\Lm]$.
Hence $\loopA$ is a $\Lm$-graded algebra with
\[\loopA^\lm = \cA^\blm \ot z^\lm\]
for $\lm\in \Lm$.
We call $\loopA$ the \textit{loop algebra} of $\cA$ relative to
the $\pi$. If we wish to emphasize the role of the grading
$\Sigma = \set{\cA^\blm}_{\blm\in\bLm}$
of $\cA$ in the loop construction, we write
$\loopA$ as $\loopASig$.
\end{definition}

\begin{remark}
\label{rem:loopgen} Let $\pi : \Lm \to \bLm$ be an epimorphism.

(i) $L_\pi$  is a functor
from the category of $\bLm$-graded algebras to the category
of $\Lm$-graded algebras.  (The morphisms in each case are the graded
homomorphisms.)

(ii) If $\cA$ is a  $\bLm$-graded algebra,
then there is a unique linear map from
$\loopA$ to  $\cA$ such that
\[u\ot z^\lm \mapsto u\]
for $\lm\in \Lm$ and $u\in\cA^\blm$. This map is an (ungraded)
algebra epimorphism of $\loopA$ onto $\cA$, and consequently the
algebra $\cA$ is a homomorphic image of the algebra $\loopA$. This
fact will be very important later in this work (see Proposition
\ref{prop:characterize}).

(iii)  It is clear that $\cA$ is a  Lie algebra
if and only if $\loopA$ is a Lie algebra.
Thus a reader whose primary interest is in Lie algebras
can chose to assume throughout this paper that all
algebras discussed are Lie algebras.  A similar remark
can be made for associative algebras, or alternative
algebras or Jordan algebras (if $k$ has characteristic
$\ne 2$).
\end{remark}

\subsection{Multiloop algebras}
\label{subsec:multi}

In this subsection
we assume that $\Lm = \Zn$ and consider multiloop algebras
graded by $\Zn$.  These are a special
case of the loop algebras just described.

If $\ell\ge 1$, an element $\zeta_\ell\in k^\times$ is called (as usual) a
\textit{primitive $\ell^\text{th}$-root of unity in $k$} if $\zeta_\ell$
generates a subgroup of $k^\times$ of order $\ell$.

\begin{definition}
\label{def:multiloop} Let
$\Lm = \Zn$.
Suppose that
$\boldm = (m_1,\dots,m_n)$ is
an $n$-tuple of positive integers
such that $k$ contains an primitive $\ell^\text{th}$
root of unity $\zeta_\ell$ (which we fix) for $\ell\in\set{m_1,\dots,m_n}$.
To construct a multiloop algebra, suppose that $\cA$ is an (ungraded) algebra
over~$k$ and
suppose that  $\sg_1,\dots,\sg_n$ is
a sequence of (pairwise) commuting finite order algebra automorphisms
of $\cA$ such that $\sg_i^{m_i} = 1$ for each $i$.
Set
\[\bLm = \bbZ/(m_1)\oplus \dots \oplus \bbZ/(m_n),\]
Then  $\cA$ has the $\bLm$-grading $\Sigma =\set{\cA^\blm}_{\blm\in\bLm}$
defined by
\begin{equation}
\label{eq:Agrading}
\cA^{(\modell_1,\dots,\modell_n)} = \set{u\in \cA \suchthat
\sg_j u = \zeta_{m_j}^{\ell_j} u \text{ for } 1\le j\le n}
\end{equation}
for $\ell_1,\dots,\ell_n\in \bbZ$, where $\modell_j = \ell_j + m_j\bbZ$ for each $j$.
We call $\Sigma$ the
\textit{$\bLm$-grading of $\cA$ determined by the automorphisms
$\sg_1,\dots,\sg_n$}.
Let
$k[z_1^{\pm 1},\dots,z_n^{\pm 1}]$ be the algebra of Laurent polynomials
over $k$ and let
\[\mloopA = \sum_{(\ell_1,\dots,\ell_n)\in \Zn}
\cA^{(\modell_1,\dots,\modell_n)}\ot z_1^{\ell_1}\cdots z_n^{\ell_n}
\subseteq \cA \otk k[z_1^{\pm 1},\dots,z_n^{\pm 1}.
\]
Then $\mloopA$ is a subalgebra of
$\cA \otk k[z_1^{\pm 1},\dots,z_n^{\pm 1}]$.
We define a $\Lm$-grading
on $\mloopA$  by setting
\begin{equation}
\label{def:multgrad}
\mloopA^{(\ell_1,\dots,\ell_n)} =
\cA^{(\modell_1,\dots,\modell_n)}\ot z_1^{\ell_1}\cdots z_n^{\ell_n}
\end{equation}
for all $(\ell_1,\dots,\ell_n)\in \Zn$.
We call the $\Zn$-graded algebra
$\mloopA$ the \textit{multiloop algebra}
of $\sg_1,\dots,\sg_n$ (based on $\cA$ and relative
to $\boldm$).

It is clear that this multiloop algebra construction
is a special case of
the general loop algebra construction described in Definition \ref{def:loop}.
Indeed, let $\pi : \Lm \to \bLm$ be the natural map defined by
\begin{equation}
\label{eq:natpi}
\pi(\ell_1,\dots,\ell_n)
= \overline{(\ell_1,\dots,\ell_n)} := (\modell_1,\dots,\modell_n),
\end{equation}
for $(\ell_1,\dots,\ell_n)\in \Lm$, and  identify
$k[z_1^{\pm1},\dots,z_n^{\pm 1}]$ with $k[\Lm]$ by means of the
map $z_1^{\ell_1}\dots z_n^{\ell_n} \mapsto z^{(\ell_1,\dots,\ell_n)}$.
Then  we have
\begin{equation}
\label{eq:mloopisaloop}
\mloopA =\loopASig.
\end{equation}
(To avoid confusion we are not abbreviating the graded algebra
$(\cA,\Sigma)$ as $\cA$ here.)
\end{definition}

\begin{remark}
\label{rem:flexibility} Our  main interest  is in multiloop
algebras.  However, as we'll see in the rest of this work, the
coordinate free point of view in the general loop construction
provides us with valuable flexibility.
\end{remark}

\begin{remark}
\label{rem:mloop}
Assume that $\Lm$, $\boldm$,  $\cA$,
$\sg_1,\dots,\sg_n$ and  $\bLm$  are as in Definition \ref{def:multiloop}.

(i) If $n=1$,
the $\bbZ$-graded algebra $M_{\boldm}(\cA,\sg_1)$
is known classically as the
loop algebra of the automorphism $\sg_1$ \cite[\S 8.2]{K}.

(ii) Although we have suppressed this from the notation (for  simplicity),
the multiloop algebra $\mloopA$ does depend on the
choice of the roots of unity $\zeta_\ell$, $\ell\in \set{m_1,\dots,m_n}$.

(iii) There is an alternate way to view the $\bLm$-grading
on $\cA$ determined by $\sg_1,\dots,\sg_n$.
To describe this, let $G=\langle
\sigma_{1},\ldots,\sigma_{n}\rangle $
and let  $\widehat{G} = \Hom(G,k^\times)$ be the
character group of $G$.
We write
\[\sg^\boldl=
{\textstyle\prod\nolimits_{i=1}^{n}} \sigma_{i}^{l_{i}}\in G\]
for  $\boldl=(l_{1},\ldots,l_{n})\in\Lm$.  Then the map
$\boldl \mapsto \sg^\boldl$
is a group epimomorphism of $\Lm$ onto $G$ that induces a
group epimorphism $\eta : \bLm \rightarrow G$ given by
\begin{equation}
\label{eq:defeta}
\eta(\boldlb) = \sg^{\boldl}
\end{equation}
for $\boldl\in\Lm$.  So $\eta$
determines a group monomorphism
$\hat{\eta} :\widehat{G}\rightarrow\widehat{\bLm}$
with
$\hat{\eta}(\chi) = \chi\circ\eta$.
Next the choice of the roots of unity
$\zeta_\ell$, $\ell\in \set{m_1,\dots,m_n}$,
defines a nondegenerate pairing
$\bLm\times \bLm\rightarrow k^{\times}$ with
\[
(\boldkb,\boldlb)\rightarrow\langle \boldkb,\boldlb\rangle =
{\textstyle\prod\nolimits_{i=1}^{n}} \zeta_{m_{i}}^{k_{i}\ell_{i}}.
\]
for  $\boldkb = (\bar k_1,\dots,\bar k_n)$ and $\boldlb = (\bar\ell_1,\dots,\bar\ell_n)$ in $\bLm$.
This pairing gives an isomorphism
$\psi:\bLm\rightarrow\widehat{\bLm}$ with $\psi(\boldlb)=\langle
\boldlb,\ \rangle $  \cite[\S 1.9, Theorem~9.2]{L}.
Consequently
\[\psi^{-1}\circ \hat{\eta} : \widehat{G} \to \bLm\]
is a group monomorphism.
But $\cA$ is naturally a
$\widehat{G}$-graded algebra (although $\widehat{G}$ is written
multiplicatively) with $\cA^{\chi}=\{a\in\cA\mid g(a)=\chi(g)a$
for all $g\in G\}$  \cite[Remark 1.3.14]{NVO2}. The monomorphism $\psi^{-1}\circ
\hat{\eta}$ transfers the $\widehat{G}$-grading of $\cA$ to a
$\bLm$-grading.  This transferred grading coincides with the $\bLm$-grading
determined by $\sg_{1},\ldots,\sg_{n}$.
Indeed, if $\boldlb = (\psi^{-1}\circ \hat{\eta})(\chi)$, where
$\chi\in \widehat{G}$, then
$\psi(\boldlb)=\chi\circ\eta$ and so
\[
\chi(\sg_{i})=(\chi\circ\eta)(\mathbf{e}_{i})=\langle
\boldlb,\mathbf{\bar{e}}_{i}\rangle =\zeta_{m_{i}}^{\ell_{i}}
\]
where $\mathbf{e}_{i}$ is the $i^\text{th}$ standard basis vector of $\Lm$. Thus,
$\cA^{\chi}=\cA^{\boldlb}$.
\end{remark}

It is sometimes convenient
to work with $\Lm$-graded algebras $\cB$ that satisfy the condition
$\langle \supp_\Lm(\cB) \rangle = \Lm$.
For example this is done in the study of various classes of tori
(see \S  \ref{sec:tori}).
Therefore, the following lemma will be useful.

\begin{lemma}
\label{lem:amplesupport}
Let $\Lm$, $\boldm$, $\cA$,
$\sg_{1},\ldots,\sg_{n}$, $\bLm$ and $\cB=\mloopA$
be as in Definition \ref{def:multiloop} and let
$\eta:\bLm\rightarrow
G=\langle \sg_{1},\ldots,\sg_{n}\rangle$
be the group epimorphism defined by
\eqref{eq:defeta}.
Then the following are
equivalent:
\begin{itemize}
    \item[(a)] $\langle \supp_{\Lm}(\cB)\rangle =\Lm$,
    \item[(b)] $\left\vert G\right\vert =m_{1}\cdots m_{n}$
    \item[(c)] $\eta$ is an isomorphism.
\end{itemize}
\end{lemma}

\begin{proof}
Clearly, (b) and (c) are equivalent. We show that
(a) and (c) are equivalent.
It is clear that
\[
\langle \supp_{\Lm}(\cB)\rangle =\Lm
\Longleftrightarrow\langle \supp_{\bar{\Lm}}(\cA)\rangle
=\bar{\Lm}\text{.}
\]
Using the $\widehat{G}$-grading of $\cA$ as in
Remark \ref{rem:mloop}(iii), let
$H=\langle \supp_{\widehat{G}}(\cA)\rangle $. Now the
natural pairing $H\times G\rightarrow k^{\times}$ is
nondegenerate. Indeed, the left kernel is trivial since $H$
consists of functions on $G$ and the right kernel is trivial since
$G$ acts faithfully on $\cA$.  Thus,
 $H = \widehat{G}$  \cite[\S 1.9, Theorem~9.2]{L}.
So $\langle \supp_{\widehat{G}}(\cA)\rangle = \widehat{G}$ and hence,
by Remark \ref{rem:mloop}(iii),
\[
\langle \supp_{\bar{\Lm}}(\cA)\rangle =\psi^{-1}
(\hat{\eta}(\widehat{G})).
\]
Therefore
\[\langle \supp_{\bar{\Lm}}(\cA)\rangle  =
\bar{\Lm}
\Longleftrightarrow
\hat{\eta}(\widehat{G})=
\widehat{\bLm}
 \Longleftrightarrow\hat{\eta}\text{ is an isomorphism}.\]
But, since $G$ and $\bLm$ are  finite,
$\hat{\eta}$ is an isomorphism if and only if $\eta$
is an isomorphism.
\end{proof}

\section{Graded-central-simple algebras}
\label{sec:graded-central-simple}

We assume again that
$\Lm$ is an arbitrary abelian group.
In  preparation for our results on the realization of
graded-central-simple algebras, we discuss in this section
some of the basic properties of these algebras.

\subsection{The centroid}
\label{subsec:centroid}

\begin{definition}
\label{def:centroid} Suppose that $\cB$ is an algebra. Let
$\Mult_k(\cB)$ be the unital subalgebra of $\End_k(\cB)$ generated by
$\set{1}\cup\set{l_a \suchthat a\in \cB}\cup\set{r_a \suchthat a\in
\cB}$, where $l_a$ (resp.~$r_a$) denotes the left (resp.~right)
multiplication operator by $a$. $\Mult_k(\cB)$ is called the
\textit{multiplication algebra} of $\cB$. Let $C_k(\cB)$ denote the
centralizer of $\Mult_k(\cB)$ in $\End_k(\cB)$.   Then $C_k(\cB)$ is a
unital subalgebra of $\End_k(\cB)$ called the  \textit{centroid}
of~$\cB$.

From now on we will usually abbreviate
$\Mult_k(\cB)$ and $C_k(\cB)$  as $\Mult(\cB)$  and $C(\cB)$ respectively.
\end{definition}

\begin{remark}
\label{rem:centre}
If $\cB$ is a unital algebra, then the map $a\mapsto l_a$ is an algebra
isomorphism of the
centre of $\cB$ onto $C(\cB)$.  (See for example \cite[\S 1]{EMO}.)
\end{remark}

Suppose that
$\cB = \lmsum\cB^\lm$
is a $\Lm$-graded algebra.
For $\lm\in \Lm$, we let
\[\End_k(\cB)^\lm = \set{e\in \End_k(\cB) \suchthat
e(\cB^\mu) \subseteq \cB^{\lm+\mu} \text{ for } \mu\in \Lm}.\]
Then $\lmsum \End_k(\cB)^\lm$
is subalgebra of $\End_k(\cB)$
that is $\Lm$-graded.
We set
\begin{equation}
\label{eq:defcentgrad}
\Mult(\cB)^\lm =
\Mult(\cB)\cap \End_k(\cB)^\lm \andd
C({\cB})^\lm  =C({\cB}) \cap \End_k({\cB})^\lm
\end{equation}
for $\lm\in \Lm$.
It is clear that
$\Mult(\cB) = \lmsum \Mult(\cB)^\lm$,
and hence $\Mult(\cB)$ is
$\Lm$-graded.  Although the centroid
is not in general $\Lm$-graded, it does have this property
in many important cases (see for
example Lemma \ref{lem:graded-simple} below).

\subsection{Graded-simplicity}
\label{subsec:graded-simplicity}

\begin{definition}  If $\cB$ is an algebra we say (as is usual)  that
$\cB$ is \textit{simple} if $\cB\cB\ne 0$ and the only ideals of $\cB$
are 0 and $\cB$. If $\cB$ is a $\Lm$-graded algebra
we say that $\cB$ is \textit{graded-simple} (or simple
as a graded algebra) if $\cB\cB\ne 0$ and the only graded ideals of $\cB$
are 0 and  $\cB$.
\end{definition}

Clearly if $\cB$ is a $\Lm$-graded algebra and $\cB\cB\ne 0$ then
\begin{equation}
\label{eq:graded-simple}
\cB \ \text{is graded-simple}\quad\iff\quad \parbox{2.4truein}{For each nonzero homogeneous element
$x\in \cB$ we have $\cB = \Mult(\cB) x$.}
\end{equation}

\begin{lemma}
\label{lem:simplecriterion}  Suppose that $\cB$ is a $\Lm$-graded algebra.  Then
\begin{equation}
\label{eq:simplecriterion}
\text{$\cB$ is simple} \iff \text{$\cB$ is graded-simple and $C(\cB)$ is a field.}
\end{equation}
Consequently, if  $C(\cB) = k1$  and $\cB$ is graded-simple then $\cB$ is
simple.
\end{lemma}

\begin{proof}  It is enough to prove \eqref{eq:simplecriterion}.
The implication ``$\Rightarrow$'' is clear (and well-known).
For the proof of ``$\Leftarrow$'', suppose that $\cB$ is graded-simple
and $C(\cB)$ is a field.  Let $\cI$ be a nonzero ideal of $\cB$.
Choose
a nonzero element
\[x = \sum_{i=1}^\ell x_i\]
in $\cI$, where $0\ne x_i\in \cB^{\lm_i}$ for all $i$ and $\lm_i\ne \lm_j$ for
$i\ne j$.  We assume that $x$ is chosen such that $\ell$ is minimum.
Since $\cB$ is graded-simple, we have
\begin{equation}
\label{eq:cyclic}
\cB = \Mult(\cB) x_1.
\end{equation}

If $\ell = 1$ then $\Mult(\cB) x = \cB$ and so $\cI = \cB$.  So we can assume that
$\ell \ge 2$.

If $m$ is a homogeneous element of $\Mult(\cB)$ and $1\le i \le \ell$,
then $mx_1 = 0 $ if and only if $mx_i = 0$ (by the minimality of $\ell$).
Consequently if $m$ and $n$ are homogeneous elements of $\Mult(\cB)$
of the same degree and $1\le i\le \ell$,
then
\begin{equation}
\label{eq:welldef}
mx_1 = nx_1 \iff mx_i = nx_i.
\end{equation}

Let $1\le i\le \ell$.  Then, by \eqref{eq:cyclic} and \eqref{eq:welldef},
there exists a well-defined map $c_i\in\End_k(\cB)^{\lm_i-\lm_1}$
such that
\[c_i(mx_1) = mx_i\]
for any homogeneous $m\in \Mult(\cB)$.
In particular
\[c_i x_1 = x_i.\]
Observe also that if $m$ and $n$
are homogeneous elements of $\Mult(\cB)$ then
$c_i(mnx_1) = mnx_i = mc_i(nx_1)$.
Hence $c_i\in C(\cB)^{\lm_i-\lm_1}$.

We now put $c = \sum_{i=1}^\ell c_i$ in which case $c \in C(\cB)$
and $x = cx_1$.  Thus, since $c$ is invertible, we have
\[\cB = c\cB = c\Mult(\cB)x_1 = \Mult(\cB) cx_1 = \Mult(\cB)x \subseteq \cI,\]
and hence $\cI = \cB$.
\end{proof}

The following is proved in \cite[Proposition 2.16]{BN}:

\begin{lemma}
Suppose that $\cB$ is a graded-simple $\Lm$-graded algebra.  Then
\label{lem:graded-simple}
\begin{itemize}
    \item[(i)] $\cB = \cB \cB$ and so $C(\cB)$ is commutative.
    \item[(ii)] $C({\cB}) = \lmsum C({\cB})^\lm$, and so $C(\cB)$ is a $\Lm$-graded  algebra.
    \item[(iii)] Each nonzero homogeneous element of $C(\cB)$ is invertible in $C(\cB)$.
    \item[(iv)] $C(\cB)^0$ is a field.
    \item[(v)] $\cB$ and $C(\cB)$ are naturally $\Lm$-graded algebras over the field $C(\cB)^0$.
    \end{itemize}
\end{lemma}

\begin{definition}
If $\cB$ is a graded-simple $\Lm$-graded algebra, we  put
\[\Gm_\Lm(\cB) := \supp_\Lm(C(\cB)) = \set{\gm\in \Lm \suchthat C(\cB)^\gm \ne 0}.\]
$\Gm_\Lm(\cB)$ is a subgroup of $\Lm$ by Lemma \ref{lem:graded-simple}(iii).
We call $\Gm_\Lm(\cB)$ the \textit{central grading group}
of $\cB$. ($\Gm_\Lm(\cB)$ is also called the centroid grading group
of $\cB$ \cite[\S6]{N}.)  From now on we will  usually abbreviate
$\Gm_\Lm(\cB)$ as $\Gm(\cB)$.
\end{definition}

\begin{remark}
\label{rem:centgroup} Suppose  that $\cB$ is a graded-simple $\Lm$-graded algebra.
If $\cB_\nu$ is obtained from $\cB$ by regrading using $\nu\in \Aut(\Lm)$
(see Definition \ref{def:regrade}),
then  $\Gm(\cB_\nu) = \nu^{-1}(\Gm(\cB))$.
\end{remark}

\subsection{Graded-centrality}
\label{subsec:graded-centrality}

\begin{definition} Suppose that $\cB$ is an algebra.
Then, $k1 \subseteq C(\cB)$, and we say that $\cB$ is \textit{central} if
$C(\cB) = k1$. We say that $\cB$ is \textit{central-simple} if
$\cB$  is central and simple. Recall that if $k$ is algebraically closed,
then  any finite dimensional simple algebra is automatically
central-simple \cite[Theorem 10.1]{J}.

Suppose next that $\cB$ is a $\Lm$-graded algebra.  Then
$k1 \subseteq C(\cB)^0 \subseteq C(\cB)$,
where $C(\cB)^0 = \set{c\in C(\cB) \suchthat
c(\cB^\lm) \subseteq \cB^{\lm} \text{ for } \lm\in \Lm}$
(see \eqref{eq:defcentgrad}).
We say that
$\cB$ is \textit{graded-central} if $C(\cB)^0 = k1$.
Further, we say that $\cB$
is \textit{graded-central-simple} if $\cB$ is graded-central
and $\cB$ is graded-simple.
\end{definition}

\begin{remark}
\label{rem:gcentral}
There  are some basic properties of graded-central-simple algebras
that suggest that they are the natural analogs of central-simple
algebras in the ungraded theory (see for example
\cite[Section 1, Chapter X]{J}).  We describe these properties here, omitting proofs
since we will not make use of the properties in this paper.

(i) If
$\cB$ is a graded-simple $\Lm$-graded $k$-algebra
and $K = C(\cB)^0$, then $K/k$ is a field extension
(see parts (iv) and (v) in Lemma \ref{lem:graded-simple}),
and $\cB$ is   naturally a graded-central-simple
$\Lm$-graded $K$-algebra.
Conversely, if
$K/k$ is a field extension and
$\cB$ is a graded-central-simple $\Lm$-graded $K$-algebra,
one can easily show that $\cB$ is
a graded-simple $\Lm$-graded $k$-algebra
and $C(\cB)^0 = K1$.

(ii) Suppose that $\cB$ is a graded-central-simple algebra over $k$
and $K/k$ is a field extension.  Then one can show that $\cB\ot_k K$ is
a graded-central-simple  $\Lm$-graded algebra  over~$K$.
\end{remark}

\begin{remark}
\label{ref:Wall}
When $\Lm = \bbZ/2\bbZ$, finite dimensional
associative unital graded-central-simple $\Lm$-graded algebras have been
classified by C.T.C.~Wall \cite{W} and they
play an important role in the theory of quadratic forms
\cite[Chapters IV and V]{Lam}.
\end{remark}

There are two cases when a graded-simple $\Lm$-graded algebra
is automatically graded-central.

\begin{lemma}
\label{lem:automaticcent}
Let $\cB$ be a graded-simple $\Lm$-graded algebra.
Suppose  either that
$\dim \cB^\lm = 1$ for some $\lm\in \Lm$ or that $k$ is algebraically closed
and $0< \dim \cB^\lm < \infty$ for some $\lm\in \Lm$.  Then
$\cB$ is graded-central-simple.
\end{lemma}

\begin{proof}  Choose $0\ne x\in \cB^\lm$.
Since $C(\cB)^0$ is a field by Lemma \ref{lem:graded-simple}(iv),
the map $c\mapsto cx$ is a linear injection of $C(\cB)^0$ into $\cB^\lm$. So
$\dim C(\cB)^0 \le \dim \cB^\lm$.  If $\dim \cB^\lm = 1$ then $C(\cB)^0 = k1$.
On the other hand if $k$ is algebraically closed and
$\dim \cB^\lm < \infty$, then $C(\cB)^0/k1$ is a finite extension and so
again $C(\cB)^0 = k1$.
\end{proof}

In view of Remark \ref{rem:gcentral}(i),
the study of graded-simple algebras over $k$ can be regarded as
equivalent to the study of graded-central-simple algebras
over extensions of $k$.  With this in mind,
we  concentrate in the rest of this paper
on the study of graded-central-simple algebras.

We first look at the structure of the centroid.  The next lemma
tells us that $C(\cB)$ is a twisted  group algebra of $\Gm(\cB)$
 over $k$ \cite[\S 1.2]{P}.

\begin{lemma} Suppose that $\cB$ is a graded-central-simple $\Lm$-graded algebra.
\label{lem:Cdecomp} Then $C(\cB)$  has a basis
$\set{c_\gm}_{\gm\in \Gm(\cB)}$ such that $c_\gm\in C(\cB)^\gm$ is a unit
of $C(\cB)$ for $\gm\in \Gm(\cB)$.  Hence
if $\gm\in \Gm(\cB)$ and $\lm\in \Lm$,
then  $\cB^{\gm+\lm}= c_\gm \cB^\lm$.
\end{lemma}

\begin{proof}  As observed in \cite[\S2.2]{BN}, this follows from Lemma
\ref{lem:graded-simple} and the fact that $C(\cB)^0 = k1$.
\end{proof}

\begin{definition}
Suppose that $\cB$ is a graded-central-simple $\Lm$-graded algebra. We say that
the centroid $C(\cB)$ of $\cB$ is \textit{split} if
\begin{equation}
\label{eq:splitdef}
C(\cB) \simeq_{\Gm(\cB)} k[\Gm(\cB)].
\end{equation}
Note that both $C(\cB)$ and $k[\Gm(\cB)]$ are $\Lm$-graded since $\Gm(\cB)$ is
a subgroup of $\Lm$.  Thus \eqref{eq:splitdef} can alternately
be written as $C(\cB) \isomLm k[\Gm(\cB)]$.  Note also that $C(\cB)$ is split
if and only  a basis $\set{c_\gm}_{\gm\in \Gm}$
for $C(\cB)$ can be chosen as in Lemma \ref{lem:Cdecomp} with the additional
property  that
\begin{equation}
\label{eq:cmult}
c_{\gm} c_{\delta}
= c_{\gm+\delta}
\end{equation}
for $\gamma,\delta\in \Gm$.
\end{definition}

If $\cB$ is an algebra, let
\[\Alg(C(\cB),k)\]
denote the set of all unital $k$-algebra homomorphisms of $C(\cB)$ into $k$.

\begin{lemma} Suppose that $\cB$ is a graded-central-simple $\Lm$-graded algebra.
\label{lem:centsplit1}
Then
\[C(\cB) \text{ is split } \iff \Alg(C(\cB),k)\ne \emptyset.\]
\end{lemma}

\begin{proof}  Since $C(\cB)$ is a twisted group ring,
this  is Exercise 17 in Chapter 1 of \cite{P}.  We include a proof
for the reader's convenience.

``$\Rightarrow$''  For this implication we can assume that
the elements $c_\gm$ in Lemma \ref{lem:Cdecomp}  satisfy \eqref{eq:cmult}.  Then the augmentation map
$c_\gm \mapsto 1$ is an element of $\Alg(C(\cB),k)$.

``$\Leftarrow$''  Suppose that $\rho\in\Alg(C(\cB),k)$.  Choose $c_\gm$, $\gm\in \Gm(\cB)$,
as in Lemma \ref{lem:Cdecomp}.  Then $\rho(c_\gm)$ is a unit in $k$
and we set
\[d_\gm = \rho(c_\gm)^{-1} c_\gm \in C(\cB)^\gm\]
for $\gm\in \Gm(\cB)$.  Now for $\gm,\delta\in \Gm(\cB)$, we have
\[d_{\gm}d_{\delta} = \rho(c_{\gm})^{-1} c_{\gm} \rho(c_{\delta})^{-1} c_{\delta}
= \rho(c_{\gm}c_{\delta})^{-1}c_{\gm}c_{\delta} =
\rho(c_{\gm+\delta})^{-1}c_{\gm+\delta}, \]
where the last equality holds since $c_{\gm}c_{\delta}$ is a nonzero
scalar multiple of $c_{\gm+\delta}$.  Hence
$d_{\gm}d_{\delta} = d_{\gm+\delta}$ as needed.
\end{proof}

There are two cases when $C(\cB)$ is always split.

\begin{lemma}
\label{lem:centsplit2}
Suppose  that
$\Lm$ is finitely generated and free, or that $k$ is algebraically closed.
If $\cB$ is a graded-central-simple $\Lm$-graded algebra,
then the centroid of
$\cB$ is split.
\end{lemma}

\begin{proof}  If $\Lm$ is finitely generated and free, then
so in $\Gm(\cB)$ and hence \eqref{eq:splitdef} is clear using
Lemma \ref{lem:Cdecomp}.  Assume next that $k$ is algebraically closed.
Then $C(\cB)$ is a commutative twisted group algebra of an abelian group
over an algebraically closed field
and so  \eqref{eq:splitdef} holds by \cite[Lemma 1.2.9(i)]{P}.
\end{proof}

\subsection{Fgc graded-central-simple algebras}
\label{subsec:fgc}

If  $\cB$ is an algebra, then $\cB$ is a (left) module over its
centroid $C(\cB)$.  We now look at this structure.

The following lemma  follows from Lemma \ref{lem:Cdecomp}
and the fact that $\Gm(\cB)$
acts freely on $\Lm$ (that is $\gm +\lm = \lm$ implies $\gamma= 0$ for $\gm\in\Gm(\cB)$
and $\lm\in \Lm$). (See  \cite[Theorem 3]{Bo} or  \cite[Lemma 2.8(ii)]{NY}).

\begin{lemma}
\label{lem:bornfree} Suppose that $\cB$ is a
graded-central-simple $\Lm$-graded algebra.
Choose a set $\Theta$
of coset representatives of $\Gm(\cB)$ in $\Lm$,
and for $\theta\in \Theta$, choose a $k$-basis $X^\theta$ for $\cB^\theta$. Using
these choices let
\begin{equation*}
X = \cup_{\theta\in \Theta} X^\theta.
\end{equation*}
Then $X$ is a homogeneous $C(\cB)$-basis for $\cB$. Hence $\cB$
is a free $C(\cB)$-module of   rank $\sum_{\theta\in \Theta}\dim_k(\cB^\theta)$
(where we interpret the sum on the right as $\infty$ if any of the terms in the sum
is infinite or if  there are infinitely many nonzero
terms in the sum).
\end{lemma}

\begin{remark}
\label{rem:submodfree}
More generally,  suppose  $\cB$ is a
graded-central-simple $\Lm$-graded algebra and
 $\cM$ is a $\Lm$-graded
$C(\cB)$-submodule of $\cB$.  Then (for the
reasons mentioned before the statement
of Lemma \ref{lem:bornfree})
$\cM$ has a homogeneous $C(\cB)$-basis and
$\cM$ is a free $C(\cB)$-module of   rank
$\sum_{\theta\in \Theta}\dim_k(\cM^\theta)$,
where
$\Theta$ is as in Lemma~\ref{lem:bornfree}.
\end{remark}

\begin{notation}
\label{not:submodfree}
Suppose that $\cB$ is a
graded-central-simple $\Lm$-graded algebra.
By the last statement in Lemma \ref{lem:Cdecomp},
we see that $\supp_\Lm(\cB)$ is the union of cosets  of  $\Gm(\cB)$ in $\Lm$.
We let
\[\supp_\Lm(\cB)/\Gm(\cB)\]
denote the set of all cosets of $\Gm(\cB)$ in $\Lm$ that are
represented by elements of $\supp_\Lm(\cB)$ (and hence consist
entirely of elements of $\supp_\Lm(\cB)$).
\end{notation}

\begin{definition}
\label{def:fgc}
If $\cB$ is an algebra we say
that $\cB$ is \textit{fgc} if
$\cB$ is finitely generated as a $C(\cB)$-module.
(The term fgc is of course an acronym
for \underline{f}initely \underline{g}enerated as a module over its \underline{c}entroid.)
\end{definition}

\begin{proposition}
\label{prop:bornmodfree}
If  $\cB$ is a
graded-central-simple $\Lm$-graded algebra, then the following are equivalent:
\begin{itemize}
    \item[(a)] $\cB$ is fgc.
    \item[(b)] $\cB$ is a free module of finite rank over $C(\cB)$.
    \item[(c)] $\supp_\Lm(\cB)/\Gm(\cB)$ is
finite and $\dim(\cB^\lm) < \infty$ for  all $\lm\in \Lm$.
\end{itemize}
Also (a), (b) and (c) are implied by
\begin{itemize}
\item[(d)] $\Lm/\Gm(\cB)$ is
finite and $\dim(\cB^\lm) < \infty$ for  all  $\lm\in \Lm$.
\end{itemize}
Moreover, if $\Lm$ is finitely generated and
$\ell\Lm \subseteq \supp_\Lm(\cB)$ for some positive
integer $\ell$, then (a), (b), (c) and (d) are equivalent.
\end{proposition}
\begin{proof}  The equivalence  of (a), (b) and (c) follow from
Lemma \ref{lem:bornfree}, and the implication ``(d) $\Rightarrow$ (c)'' is trivial.
It remains to show that (c)
implies (d) when the additional assumptions on $\Lm$ and $\supp_\Lm(\cB)$
hold.  Let $\bLm = \Lm/\Gm(\cB)$ and let $\bar{\phantom a} : \Lm \to \bLm$ be the
canonical projection.
Set $T = \supp_\Lm(\cB)$.  Then, since (c) holds, $\bar T$ is finite.
But $\ell \bLm \subseteq \bar T$ and so $\ell \bLm$ is finite.
On the other hand, since $\Lm$ is finitely generated, $\bLm/\ell\bLm$ is finite.  Therefore
$\bLm$ is finite and (d)  holds.
\end{proof}

\begin{remark}  Suppose that $\Lm$ is finitely generated.
The additional assumption that $\ell\Lm \subseteq \supp_\Lm(\cB)$
for some positive integer $\ell$ holds for many important classes of graded-central-simple
$\Lm$-graded algebras such as centreless Lie tori,
associative tori, alternative tori
and Jordan tori (see Section \ref{sec:tori}).
Thus for any $\Lm$-graded algebra $\cB$
in one of these classes, (a) is equivalent to (d) in
Proposition  \ref{prop:bornmodfree}.
\end{remark}

\section{Graded simplicity and centrality for loop algebras}
\label{sec:gradedloop}
Suppose in  this section that $\Gm$ is a subgroup of
an arbitrary abelian group $\Lm$.
Suppose further that $\pi : \Lm \to \bLm$ is a (group) epimorphism
of $\Lm$ onto an abelian group $\bLm$ such that
\[\kerpi = \Gm.\]

In this section we investigate centrality and simplicity
of the loop algebra $\loopA$.

\subsection{Preliminary lemmas}
\label{subsec:loopproperties}
\begin{lemma}
\label{lem:loopgraded-simple} Suppose that $\cA$ is a
$\bLm$-graded algebra.  Then
\[\textit{$\cA$ is graded-simple} \iff \textit{$\loopA$ is
graded-simple.}\]
\end{lemma}

\begin{proof} ``$\Leftarrow$''  If $\cA$ is not graded-simple, then $\cA$ has a nonzero
proper graded ideal $\cI$.  In that case, $L_\pi(\cI)$ is
a nonzero proper graded ideal of $\loopA$, and so $\loopA$ is not
graded-simple.

``$\Rightarrow$''  Suppose that $\cA$ is graded-simple.
Note first that the set $\set{u\in \cA \suchthat u\cA + \cA u = 0}$
is a proper graded ideal of $\cA$ and hence it is 0.  Thus, if
$0\ne u\in \cA$ then $u\cA + \cA u \ne 0$.

Let $\cI$ be a nonzero graded ideal of $\loopA$.
For $\lm\in \Lm$ we let
\[\cS^{\blm} := \set{u\in \cA^{\blm} \suchthat u\ot z^{\gm+\lm}\in \cI
\text{ for all } \gm\in\Gm },\]
and we set $\cS = \sum_{\blm\in\bLm}\cS^{\blm}$. Then $\cS$ is a graded ideal
of $\cA$.

Now choose
$0\ne u_1\ot z^\mu
\in \cI$, where $\mu\in \Lm$ and $u_1\in \cA^{\bar \mu}$.
Then by the note at the beginning of this proof,
we may choose $\lm\in\Lm$ and $u_2\in \cA^{\blm}$ such that
$u_1u_2 \ne 0$ or $u_2u_1\ne 0$.  We assume that $u_2u_1 \ne 0$ (the other case
being similar).  For any $\gm\in \Gm$, the element
$(u_2\ot z^{\gm+\lm})(u_1\ot z^\mu) = u_2u_1 \ot z^{\gm+\lm+\mu}$
is in $\cI$.
Hence $u_2u_1\in\cS^{\blm + \bmu}$ and so $\cS\ne 0$.  Consequently $\cS = \cA$.
Thus for all $\lm\in \Lm$ and all $u\in \cA^\blm$, we have
$u\ot z^\lm\in\cI$.  Therefore $\cI = \loopA$.
\end{proof}

Suppose that  $\cA$ is a graded-simple $\bLm$-graded algebra.
By Lemma \ref{lem:graded-simple}(ii), $C(\cA)$
is a $\bLm$-graded algebra and so $L_\pi(C(\cA))$ is a
$\Lm$-graded algebra. On the other hand, by Lemma \ref{lem:loopgraded-simple},
$\loopA$ is a graded-simple $\Lm$-graded algebra and so
$C(\loopA)$ is a $\Lm$-graded algebra. We let
\[\psi : L_\pi(C(\cA)) \to C(\loopA)\]
be the unique $k$-linear map
such that
\begin{equation}
\label{eq:defpsi}
\big(\psi(c\ot z^\lm)\big)(u\ot z^\mu) = c(u) \ot z^{\lm+\mu}
\end{equation}
for $\lm,\mu\in\Lm$, $c\in C(\cA)^\blm$, $u\in \cA^\bmu$.
It is easy to check that $\psi$ is a homomorphism of $\Lm$-graded
algebras.

The  first part of the following  lemma (with  weaker assumptions) was proved for
classical loop algebras in \cite[Proposition 4.11]{ABP}.

\begin{lemma}
\label{lem:psiiso}
Suppose that
$\cA$ is a graded-simple $\bLm$-graded algebra.
Then the map $\psi : L_\pi(C(\cA)) \to C(\loopA)$ defined by
\eqref{eq:defpsi} is an isomorphism of $\Lm$-graded algebras.
Moreover,
\begin{equation}
\label{eq:Gmloop}
\Gm_\Lm(\loopA) = \set{\lm\in \Lm \suchthat \blm \in \Gm_\bLm(\cA)}.
\end{equation}
\end{lemma}

\begin{proof} $\ker(\psi)$ is a graded ideal of $L_\pi(C(\cA))$ and it is clear
that $\ker(\psi)^\lm= 0$ for $\lm\in \Lm$.  Thus $\psi$ is a monomorphism.

To see that $\psi$ is onto, let $d\in C(\loopA)^\lm$ where $\lm\in \Lm$.
Then for $\mu\in \Lm$, there exists a unique map
\[c_\mu : \cA^\bmu \to \cA^{\bmu+\blm}\]
such that
\begin{equation}
\label{eq:dcomp}
d(u \ot z^\mu) = c_\mu(u)\ot z^{\mu+\lm}
\end{equation}
for $u\in \cA^\bmu$.  Then for $\mu_1,\mu_2\in\Lm$, $u_1\in \cA^{\bmu_1}$,
$u_2\in \cA^{\bmu_2}$, we have
\begin{align*}
c_{\mu_1+\mu_2}(u_1u_2)\ot z^{\mu_1+\mu_2+\lm} &= d(u_1u_2\ot z^{\mu_1+\mu_2})
= (u_1\ot z^{\mu_1})d(u_2\ot z^{\mu_2})\\
&= (u_1\ot z^{\mu_1})(c_{\mu_2}(u_2)\ot z^{\mu_2+\lm})
= u_1c_{\mu_2}(u_2)\ot z^{\mu_1+\mu_2+\lm}.
\end{align*}
Hence
\begin{equation}
\label{eq:cmu}
c_{\mu_1+\mu_2}(u_1u_2) = u_1c_{\mu_2}(u_2) \quad\text{ and similarly }\quad
c_{\mu_1+\mu_2}(u_1u_2) = c_{\mu_1}(u_1)u_2
\end{equation}
for $\mu_1,\mu_2\in\Lm$, $u_1\in \cA^{\bmu_1}$,
$u_2\in \cA^{\bmu_2}$.
Thus if $\mu_1,\mu_2\in\Lm$, $\gm\in \Gm$, $u_1\in \cA^{\bmu_1}$,
$u_2\in \cA^{\bmu_2}$, we have
\[c_{\mu_1+\mu_2+\gm}(u_1u_2) = c_{(\mu_1+\gm)+\mu_2}(u_1u_2) = u_1c_{\mu_2}(u_2)
= c_{\mu_1+\mu_2}(u_1u_2).
\]
Consequently if $\mu\in \Lm$ and $\gm\in \Gm$, the maps
$c_{\mu+\gm}$ and $c_\mu$
agree on all elements of $(\cA\cA)\cap \cA^\bmu$.  Since $(\cA\cA)\cap \cA^\bmu = \cA^\bmu$,
it follows that
$c_{\mu+\gm} = c_\mu$
for $\mu\in \Lm$, $\gm\in \Gm$.
This equation tells us that there is a unique well-defined map $c\in\End_k(\cA)$
such that
$c(u) = c_\mu(u)$
for  $u\in \cA^\bmu$ and  $\mu\in \Lm$.
Clearly $c\in\End_k(\cA)^\blm$ and then by \eqref{eq:cmu}, we have
$c\in C(\cA)^\blm$.  Hence $c\ot z^\lm \in L_\pi(C(\cA))^\lm$.
Finally, by \eqref{eq:dcomp}, we have
$d(u\ot z^\mu) = c(u) \ot z^{\mu+\lm}$
for $\mu\in \Lm$ and $u\in \cA^\bmu$.  Thus $\psi(c\ot z^\lm) = d$.

So we have proved the first statement.  Hence we have
$C(\loopA) \isomLm L_\pi(C(\cA)) = \sum_{\lm\in\Lm}C(\cA)^\blm \ot z^\lm$,
which implies \eqref{eq:Gmloop}.
\end{proof}

\begin{lemma}
\label{lem:Bcentral} Suppose that  $\cA$ is a $\bLm$-graded
algebra.
\begin{itemize}
    \item[(i)] If $\cA$ is graded-simple, then
    \[\textit{$\cA$ is graded-central} \iff \textit{$\loopA$ is graded-central}.\]
    \item[(ii)] If $\cA$ is graded-central-simple, then
    \[\textit{$\cA$ is central-simple} \iff \textit{$\Gm(\loopA)=\Gm$}.\]
    \item[(iii)] If $\cA$ is central-simple, then
\begin{equation}
\label{eq:CloopA}
    C(\loopA) = \spann_k\set{l_{1\ot z^\gm} \suchthat \gm\in \Gm}\isomLm k[\Gm],
\end{equation}
    where $l_{1\ot z^\gm}$ denotes left multiplication by $1\ot z^\gm$.  In particular,
    the centroid of $\loopA$ is split
    \end{itemize}
\end{lemma}

\begin{proof}
(i):  Suppose that $\cA$ is graded-simple. Then $L_\pi(C(\cA))^0 = C(\cA)^{\bar 0}\ot z^0$. So
$L_\pi(C(\cA))^0$ and $C(\cA)^{\bar 0}$ are isomorphic as unital $k$-algebras.
By Lemma \ref{lem:psiiso}, $C(\loopA)^0$ and
$C(\cA)^{\bar 0}$ are thus isomorphic as unital $k$-algebras and (i) follows from this.

(ii):  Suppose that $\cA$ is graded-central-simple.
Then,
\begin{align*}
\Gm(\loopA) = \Gm &\iff \Gm_\bLm(\cA) = \set{\bar 0} &&\text{(by \eqref{eq:Gmloop})}\\
&\iff C(\cA) = C(\cA)^{\bar 0}\\
&\iff C(\cA) = k1 && \text{(since $\cA$ is graded-central)}\\
&\iff \cA \text{ is central}\\
&\iff \cA \text{ is central-simple,}
\end{align*}
where the last equivalence follows from the second statement in
Lemma \ref{lem:simplecriterion}.

(iii): Observe
that
\[\textstyle
L_\pi(C(\cA))\\
= \sum_{\lm\in\Lm}C(\cA)^\blm \ot z^\lm
= \sum_{\gm\in\Gm} k \ot z^\gm
=\spann_k\set{1\ot z^\gm \suchthat \gm\in \Gm}.
\]
Hence we have \eqref{eq:CloopA} by Lemma \ref{lem:psiiso}.
\end{proof}

\subsection{The classes $\classA$ and $\classB$}
\label{subsec:classes}

For use
here and in the rest of the
paper we now introduce a  class $\classA$ of $\bLm$-graded algebras
and a class $\classB$ of  $\Lm$-graded algebras.
As the notation suggests, the class $\classA$ depends on just the group $\bLm$,
while the class $\classB$ depends both on the group $\Lm$
and the subgroup~$\Gm$.

\begin{definition}
\label{def:class}
(i) We let $\classA$ be the class of $\bLm$-graded algebras $\cA$
such that  $\cA$ is central-simple as an algebra.

(ii) We let
$\classB$ be the class of $\Lm$-graded algebras $\cB$
 such that $\cB$ is graded-central-simple, the centroid of $\cB$ is split
and $\Gm(\cB) = \Gm$.  Equivalently
$\classB$ is the class of $\Lm$-graded algebras $\cB$
such that $\cB$ is graded-central-simple and $C(\cB) \isomLm  k[\Gm]$.
\end{definition}

\begin{remark}
\label{rem:class}
It is clear that the class $\classA$ is closed under graded-isomorphism
That is, if $\cA$ and $\cA'$ are $\bLm$-graded algebras such that
$\cA\in \classA$ and $\cA\isombLm\cA'$, then we also have $\cA'\in \classA$.
Similarly $\classB$ is closed under graded-isomorphism.
\end{remark}

In the next proposition we use the loop construction and
the previous lemmas to establish a relationship between the classes $\classA$ and $\classB$
using the loop construction.  This relationship will be explored in
more detail in \S \ref{sec:correspondence}.

\begin{proposition}
 \label{prop:loopGCS}
Let $\cA$ be a
$\bLm$-graded algebra.  Then the following statements are equivalent:
\begin{itemize}
    \item [(a)] $\cA\in \classA$.
    \item [(b)] $\loopA \in \classB$.
    \item [(c)] $\loopA$ is a graded-central-simple with central grading group $\Gm$.
\end{itemize}
\end{proposition}
\begin{proof}
``(a) $\Rightarrow$ (b)'' follows from Lemmas \ref{lem:loopgraded-simple} and
Lemmas \ref{lem:Bcentral},  whereas ``(b) $\Rightarrow$ (c)'' is
trivial.  Finally ``(c) $\Rightarrow$ (a)'' follows by Lemmas
\ref{lem:loopgraded-simple} and  Lemmas~\ref{lem:Bcentral}.
\end{proof}

\section{Central specializations and central images}
\label{sec:central-spec}
We have seen  in Section \ref{sec:loop} how to pass from a
$\bLm$-graded algebra to a $\Lm$-graded algebra  using the loop
construction. In order to provide an inverse for this construction
(in a sense to be made precise),  we study in this section
certain algebra homomorphisms, called central specializations,
from $\Lm$-graded algebras onto $\bLm$-graded algebras.

Throughout  the section we assume again that  $\Gm$ is a subgroup
of
an arbitrary abelian group $\Lm$
and that $\pi : \Lm \to \bLm$ is a group epimorphism such  that $\Gm
= \kerpi.$

\subsection{Definitions}
\label{subsec:defnspec}

\begin{definition}
\label{def:spec} Let
$\cB$ be a $\Lm$-graded algebra and  let
$\rho\in \algC$. A \textit{$\rho$-specialization}
of $\cB$ is a nonzero algebra epimorphism
$\ph : \cB \to \cA$ onto a $\bLm$-graded algebra $\cA$ such that the
following two conditions hold:
\begin{itemize}
\item[(a)] $\ph(\cB^\lm) \subseteq \cA^\blm$ for $\lm\in \Lm$.
\item[(b)] $\ph(cx) = \rho(c)\ph(x)$ for $c\in C(\cB)$, $x\in \cB$.
\end{itemize}
If $\ph : \cB \to \cA$ is a $\rho$-specialization, we call
$\cA$ a \textit{$\rho$-image of $\cB$}.

A \textit{central specialization} of $\cB$
is a map
$\ph : \cB \to \cA$ that is a $\rho$-specialization
of $\cB$ for some $\rho\in \algC$.
Similarly a \textit{central image} of $\cB$ is
a $\bLm$-graded algebra that is
a $\rho$-image of $\cB$ for some $\rho\in \algC$.

Of course all of these definitions are made relative to the fixed
epimorphism $\pi : \Lm \to \bLm$.
\end{definition}

\begin{remark}
\label{rem:histspec}  Suppose  that
$\cA$ is a finite dimensional simple Lie algebra
over an algebraically closed field of characteristic 0.
Although  not formulated as in Definition
\ref{def:spec}, central specializations of classical loop algebras
of $\cA$ were used by  V.~Kac  in the classification of automorphisms of
finite order of $\cA$  \cite[Theorem~8.6]{K}.
\end{remark}

\begin{remark}
\label{rem:inducegrading}
Let
$\cB$  be a $\Lm$-graded algebra
and let
$\rho\in \algC$.

(i)
Since $\cB$ is $\Lm$-graded, $\cB$ has a natural $\bLm$-grading
defined by
\begin{equation}
\label{eq:gradquot}
\textstyle
\cB^\blm = \sum_{\mu\in \Lm,\ \bmu = \blm}
\cB^\mu
\end{equation}
for $\lm\in \Lm$. (See \cite[p.~3]{NVO2}.) Then assumption (a) in
Definition \ref{def:spec} says that $\ph$ is a $\bLm$-graded map.

(ii)  Suppose $\cB'$ is another $\Lm$-graded algebra
and $\beta : \cB' \to \cB$ is an isomorphism of $\Lm$-graded algebras.
Then $\beta$ induces an algebra isomorphism
\[C(\beta) : C(\cB') \to C(\cB)\]
defined by $C(\beta)(c') = \beta\circ c' \circ\beta^{-1}$ for $c'\in C(\cB')$.
It follows that $C(\cB') = \bigoplus_{\gm\in\Gm} C(\cB')^\gm$
and that $C(\beta)$ is an isomorphism of $\Lm$-graded algebras.
Moreover, if $\ph : \cB \to \cA$ is a $\rho$-specialization
of $\cB$, then $\ph\circ\beta$ is a $\rho\circ C(\beta)$-specialization
of $\cB'$.  Consequently if $\cA$ is a central image of $\cB$
then $\cA$ is also a central image of $\cB'$.

(iii)  If $\cA$ is a $\rho$-image of $\cB$ and
$\cA'$ is $\bLm$-graded algebra such that $\cA \isombLm \cA'$,
then $\cA'$ is also a $\rho$-image of $\cA$.
\end{remark}

\begin{example}
\label{ex:universal}
Let   $\cB$ be a $\Lm$-graded algebra
which satisfies the following conditions:
\begin{itemize}
\item[(a)] $C(\cB)$ is commutative and
$C(\cB) = \bigoplus_{\gm\in\Gm} C(\cB)^\gm$
(where $C(\cB)^\gm$ is defined by \eqref{eq:defcentgrad}
for $\gm\in \Gm$).
\item[(b)] $\cB$ is a nonzero free $C(\cB)$-module (under the natural action).
\end{itemize}
(Note that by Lemmas
\ref{lem:graded-simple}(ii) and
\ref{lem:bornfree},
these conditions are satisfied if $\cB\in \classB$.)
Suppose
$\rho\in \algC$.  Let
\[\kerrhoB = \spann_k\{cx \suchthat c\in \ker(\rho),\ x\in \cB\}.\]
Then $\kerrhoB$ is an ideal of $\cB$ (as an algebra).
Also, regarding $\cB$ as $\bLm$-graded as in Remark \ref{rem:inducegrading}(i),
we have, using assumption (a), that
\[\ker(\rho)\cB^\blm \subseteq C(\cB) \cB^\blm
\subseteq
\textstyle
\big(\sum_{\gm\in\Gm} C(\cB)^\gm\big)
\big(\sum_{\mu\in \Lm,\ \bmu = \blm}\cB^\mu\big) \subseteq \cB^\blm
\]
for $\blm\in \bLm$.  It follows from this that
$\kerrhoB$ is a $\bLm$-graded ideal of $\cB$.  Thus the quotient
algebra
\[\unim\]
has the natural structure of a $\bLm$-graded algebra.
Observe also that $\ker(\rho) \ne C(\cB)$ and so,
by assumption (b), $\kerrhoB \ne \cB$.  (Actually it would be
enough to assume in place of (b) that $\cB$ is a faithfully flat
$C(\cB)$-module \cite[\S I.3.1]{Bour}.)  Thus $\unim \ne 0$.
Finally, let $\projrho : \cB \to \unim$ be  the canonical projection defined
by
\[\projrho(x) = x + \kerrhoB\]
for $x\in \cB$.
Note that since
$(c-\rho(c)1)x\in \ker(\rho)\cB$,
we have
\[\projrho(cx) = \rho(c)\projrho(x)\]
for $c\in C(\cB)$, $x\in \cB$.
Thus $\projrho$ is a $\rho$-specialization of $\cB$.
We call  $p_\rho$ the \textit{universal
$\rho$-specialization of $\cB$}.  This terminology is
justified by the following  lemma.
\end{example}

\begin{lemma}
\label{lem:universal}
Let $\cB$ be a $\Lm$-graded
algebra satisfying assumptions (a) and (b) of
Example \ref{ex:universal} and let
$\rho\in \algC$.
If $\ph : \cB \to \cA$ is an arbitrary $\rho$-specialization
of $\cB$, then there is a unique $\bLm$-graded epimorphism
$\kappa : \unim \to \cA$ such that $\ph = \kappa \circ \projrho$.
\end{lemma}
\begin{proof} Now  $\ker(\rho)\cB \subseteq \ker(\ph)$.
So there is an induced algebra homomorphism $\kappa :
\unim \to \cA$ such that
\[\kappa(x + \ker(\rho)\cB) = \ph(x)\]
for $x\in \cB$.  It is clear that $\kappa$ is $\bLm$-graded,
and since $\ph$ is surjective, $\kappa$ is surjective.
Also, $\ph = \kappa \circ \projrho$ by definition of $\kappa$.
Finally, the uniqueness of $\kappa$ is clear since $\projrho$
is surjective.
\end{proof}

\begin{remark}
\label{rem:universal} Let
$\cB$ be a $\Lm$-graded
algebra satisfying assumptions (a) and (b) of
Example \ref{ex:universal} and let
$\rho\in \algC$.
Then the universal $\rho$-specialization
of $\cB$
has an alternate interpretation that was suggested to us by
Ottmar Loos.  Indeed, we may regard $k$ as an algebra $k_\rho$
over $C(\cB)$ by means of
the action $(c,a)\mapsto \rho(c)a$ for $c\in C(\cB)$ and $a\in k$.
Then
\[\unten\]
is an algebra over $k$.  Furthermore, as we have seen, it follows from assumption (a) that
$\cB = \bigoplus_{\blm\in\bLm}\cB^{\blm}$ is a decomposition of $\cB$
as the  direct sum of $C(\cB)$-modules. Hence $\unten$ is a $\bLm$-graded algebra with
$(\unten)^\blm := \set{x\ot 1 \suchthat x\in \cB^\blm}$ for $\blm\in \bLm$.
Moreover, $\unten$ is nonzero
by assumption (b). Finally one checks easily that the map $x +
\ker(\rho)\cB \mapsto x\ot 1$ is a $\bLm$-graded algebra isomorphism
of $\unim$ onto $\unten$.  If we regard this map as an
identification, then the universal $\rho$-specialization $p_\rho :
\cB \to \unten$ is given by $p_\rho(x) = x\ot 1$ for   $x\in
\cB$.
\end{remark}

\subsection{Central specializations and images of algebras in $\classB$}
\label{subsec:propspec}

If $\cB$  is in $\classB$, then by
Lemma \ref{lem:centsplit1}
we have $\algC \ne \emptyset$.  Moreover,
$\cB$ satisfies assumptions (a) and (b) of Example
\ref{ex:universal}, and so, for $\rho\in \algC$,
we can construct
the universal $\rho$-specialization
$\projrho: \cB \to \unim$ of $\cB$.
In part (i) of the next proposition,
we see that $\projrho$ is
unique
$\rho$-specialization  of $\cB$.

\begin{proposition}
\label{prop:spec} Suppose that $\cB\in\classB$, $\rho\in \algC$, $\cA$
is a $\bLm$-graded algebra and
$\ph : \cB \to \cA$ is a $\rho$-specialization of $\cB$.
Then
\begin{itemize}
    \item[(i)] There exists a unique $\bLm$-graded isomorphism
    $\kappa : \unim \to \cA$ such that $\ph = \kappa \circ \projrho$.
     \item[(ii)]  If $X$ is a homogeneous $C(\cB)$-basis for $\cB$ chosen as
    in Lemma \ref{lem:bornfree}, then $\ph$ maps $X$ bijectively onto a $k$-basis
    $\ph(X)$ of $\cA$.
    \item[(iii)]  For $\lm\in \Lm$,  $\ph$ restricts to a linear bijection
    of $\cB^\lm$ onto $\cA^\blm$.
    \item[(iv)] $\loopA \isomLm \cB$.
    \item[(v)] $\cA\in \classA$.
    \end{itemize}
\end{proposition}

\begin{proof}
We first prove statements (ii)-(v) for the universal $\rho$-specialization.
So in this part of the proof
we assume that $\cA = \unten$ and
$\ph : \cB \to
\unten$ is given by $\ph(x) = x\ot 1$
(see Remark \ref{rem:universal}).

(ii): This is a general property of the tensor product $\unten$
(see \cite[Proposition 4.1]{L}).

(iii):  This follows from (ii) since we can choose the
$C(\cB)$-basis $X$ for $\cB$ so that $X$ contains a $k$-basis for  $\cB^\lm$.

(iv): Define $\omega : \cB \to \loopA$ by
\[\omega(x) = \ph(x)\ot z^\lm\]
for $x\in \cB^\lm$ and  $\lm\in\Lm$.  Then $\omega$ is a nonzero homomorphism
of $\Lm$-graded algebras.  Since $\cB$ is graded-simple,
$\omega$ is a monomorphism. Finally, by (iii),
we have $\omega(\cB^\lm) = \cA^\blm \ot z^\lm$ for  $\lm\in\Lm$, and so
$\omega$ is surjective.

(v):  By (iv), we have $\loopA \isomLm \cB$.
So $\loopA$ is graded-central-simple with central grading group
$\Gm$.
Thus, by Proposition  \ref{prop:loopGCS},  $\cA\in \classA$.

To complete the proof of the proposition, we now
assume that $\ph : \cB \to \cA$ is an arbitrary
$\rho$-specialization.  Since we now know
that $\unim$ is graded-simple (in fact it is  central-simple by (v) in the universal case),
(i) follows from the
last statement of Lemma \ref{lem:universal}.
(ii)--(v) then follow from the corresponding statements
for the universal $\rho$-specialization.
\end{proof}

By Proposition \ref{prop:spec}(i)
we have the following:

\begin{corollary}
\label{cor:uniquecent}
Suppose $\cB\in\classB$ and
$\rho\in \algC$.  Then $\unim$ is the unique $\rho$-image of $\cB$ up
to $\bLm$-graded isomorphism.
\end{corollary}

\subsection{Similarity for $\bLm$-graded algebras}
\label{subsec:similarity}
In the  subsection after this we will look
at central images  of $\cB \in \classB$ corresponding
to different homomorphisms in $\algC$.
In preparation for this we look first at a notion of
similarity (relative to $\pi$) for $\bLm$-graded
algebras.

Since $\pi : \Lm\to \bLm$ is surjective,
$\pi$ has a right inverse as a map of sets.
For the  rest of \S 6,
\textit{we fix a choice $\xi$ of such a right inverse}.
So $\xi : \bLm \to \Lm$ is a map of sets such that
\[\pi\circ\xi = 1_{\bLm}.\]

\begin{definition}
\label{def:twist}
Let $\chi$ be a character of $\Gamma$, i.e. $\chi\in\charG$.
Also let $\cA$ be a $\bLm$-graded algebra. We define
a $\bLm$-graded algebra $\Achi$ as follows.  As a $\bLm$-graded
vector space $\Achi = \cA$.  Further, the product $\dotchi$ on
$\cA_\chi$ is defined by
\[u\dotchi v = \chi(\xi(\blm) + \xi(\bmu) - \xi(\blm + \bmu)) uv\]
for $\blm,\bmu\in \bLm$, $u\in \cA^\blm$, $v\in \cA^\bmu$.
We call $\Achi$ the \textit{twist of $\cA$ by $\chi$}.
\end{definition}

\begin{remark}
\label{rem:twistind}
Suppose that $\chi$ and $\cA$ are as in Definition
\ref{def:twist}.
It is easy to check that, up to $\bLm$-graded isomorphism,
$\Achi$ is independent of the choice of the right inverse $\xi$ for $\pi$.
\end{remark}

\begin{remark}
\label{rem:assocLietwist}  Suppose that $\cA$ is a $\bLm$-graded algebra.
If $\cA$ is a Lie or associative algebra, then one easily checks directly that any twist of
$\cA$ is a Lie algebra or associative algebra respectively.   Note also that
if $\cA$ is unital then any twist $\cA_\chi$
of $\cA_\chi$ is unital.  Indeed, up to
a $\bLm$-isomorphism of $\cA_\chi$, we can choose the right inverse
$\xi$ with $\xi(\bar 0) = 0$.  In this case,
the identity element $1$ of $\cA$ has $1\in\cA^{\bar 0}$ and
$1$ is also an identity of $\cA_\chi$.
\end{remark}

Twists of $\cA$ have the following properties:

\begin{lemma}
\label{lem:twistprop}  Suppose that $\cA$ is a $\bLm$-graded algebra.
\begin{itemize}
\item[(i)] If $\cA'$ is a $\bLm$-graded algebra
such that $\cA\isombLm \cA'$, then
$\cA_\chi \isombLm \cA'_\chi$ for $\chi\in\charG$.
\item[(ii)]  If $1\in\charG$ is the trivial character (that is
$1(\gm) = 1$ for all $\gm\in\Gm$),
then $\cA_1 = \cA$.
\item[(iii)] If $\chi_1,\chi_2\in \charG$, then
$(\cA_{\chi_1})_{\chi_2} = \cA_{\chi_1\chi_2}$.
\item[(iv)]  If $\chi\in\charG$ extends to a character
of $\Lm$, then $\Achi \simeq_\bLm \cA$.
\item[(v)]  If $k$ is
algebraically  closed, then
$\Achi \simeq_\bLm \cA$
for any $\chi\in \charG$.
\end{itemize}
\end{lemma}

\begin{proof}  (i) and (ii) are clear, and (iii) is easily checked.

For (iv), suppose that $\chi\in \charG$ and $\chi$ extends
to a character $\psi$ of $\Lm$.
Then for $\blm,\bmu\in \bLm$, $u\in \cA^\blm$, $v\in \cA^\bmu$,
we have
\[u\dotchi v = \psi(\xi(\blm))\psi(\xi(\bmu))\psi(\xi(\blm+\mu))^{-1}  uv.\]
Hence the map
defined by
$u \mapsto \psi(\xi(\blm))^{-1}u$ for $\blm\in\bLm$, $u\in \cA^\blm$
is a $\bLm$-graded isomorphism of $\cA$
onto $\Achi$.

For (v), suppose that $k$ is algebraically closed.
Then by \cite[Lemma 1.2.7]{P},
any character of $\Gm$  extends to a character of $\Lm$.  So (v)
follows from (iv).
\end{proof}

\begin{remark}
\label{rem:cohom} Twists of $\cA$ and the preceding lemma have a
cohomological interpretation using the exact  sequence
$H^1(\Lm,k^\times) \to  H^1(\Gm,k^\times) \to H^2(\bLm,k^\times)$
arising from the exact sequence $0\to \Gm \to \Lm \to \bLm \to 0$
(with trivial actions on $k^\times$) \cite[\S 2.6]{S}.
Since we will not make use of this here, we omit the details.
\end{remark}

\begin{definition}
\label{def:similar}
If $\cA$ and $\cA'$ are $\bLm$-graded algebras, we say that
$\cA$ and $\cA'$ are \textit{similar relative to $\pi$}, written
$\cA \sim_\pi \cA'$, if
$\cA' \isombLm \cA_\chi$
for some
$\chi\in \charG$.
\end{definition}

\begin{remark}
\label{rem:twist}
(i)
The relation $\sim_\pi$ depends on the group epimorphism
$\pi : \Lm \to \bLm$ (with kernel $\Gamma$) but not on
the choice of the right inverse
$\xi$ for $\pi$.  (See Remark \ref{rem:twistind}.)

(ii)  It follows from
parts (i)--(iii) of Lemma \ref{lem:twistprop}  that $\sim_\pi$ is an equivalence
relation on the class of $\bLm$-graded algebras.

(iii)  Suppose that $\cA$ and $\cA'$ are $\bLm$-graded algebras. By
Lemma \ref{lem:twistprop}(ii) we see that
\[\cA \isombLm \cA' \Rightarrow \cA \sim_\pi \cA'.\]
Moreover, if $k$ is algebraically closed, then   by Lemma \ref{lem:twistprop}~(v), we have
\[\cA \isombLm \cA'  \iff \cA \sim_\pi  \cA'.\]
\end{remark}

\subsection{Comparing central images}
\label{subsec:compare}

We now look at central images  of $\cB \in \classB$ corresponding
to different homomorphisms in $\algC$.

\begin{lemma}
\label{lem:rhochi} Suppose  that $\cB\in\classB$ and $\rho\in \Alg(C(\cB),k)$.
\begin{itemize}
    \item[(i)] For $\chi\in\Hom(\Gamma,k^{\times})$, set $\rho_\chi(c)=\chi(\gamma)\rho(c)$
for $\gamma\in\Gamma$ and $c\in C(\cB)^{\gamma}$. Then $\chi\rightarrow\rho_\chi$
    is a bijection of $\Hom(\Gamma,k^{\times})$ onto $\Alg(C(\cB),k)$.

\item[(ii)] If $\chi\in\Hom(\Gamma,k^{\times})$ and $\cA$ is a $\rho$-image of $\cB$, then
$\cA_{\chi}$ is a $\rho_\chi$-image of $\cB$.

\item[(iii)]If $\chi\in\Hom(\Gamma,k^{\times})$, $\cA$ is a $\rho$-image of $\cB$ and
$\cA^{\prime}$ is a $\rho_\chi$-image of $\cB$, then
$\cA^{\prime}\simeq_{\bar{\Lm}}\cA_{\chi}$.
\end{itemize}
\end{lemma}

\begin{proof} (i) is clear since
$C(\cB)$ is isomorphic to $k[\Gm]$ by Lemma \ref{lem:centsplit1}.

(ii): Suppose $\chi\in\Hom(\Gamma,k^{\times})$ and $\cA$ is a $\rho$-image of $\cB$.
Then we have a $\rho$-specialization $\ph : \cB \to \cA$.  Define
$\ph' : \cB \to \cA_\chi$ by
\[\ph'(x) = \chi(\lm-\xi(\blm))\ph(x)\]
for $\lm\in \Lm$, $x\in \cB^\lm$.
Then $\ph'$ is nonzero, surjective and $\bLm$-graded.  Also
if $\lm,\mu\in \Lm$, $x\in \cB^\lm$ and $y\in \cB^\mu$, we have
\begin{align*}
\ph'(xy) &= \chi(\lm+\mu - \xi(\blm+\bmu))\ \ph(xy)\\
&= \chi(\lm+\mu - \xi(\blm+\bmu))\ \ph(x)\ph(y)\\
&= \frac {\chi(\lm+\mu - \xi(\blm+\bmu))}
{\chi(\xi(\blm) + \xi(\bmu) - \xi(\blm + \bmu))}\
\ph(x)\dotchi\ph(y)\\
&= \chi(\lm+\mu - \xi(\blm) - \xi(\bmu))\ \ph(x)\dotchi\ph(y)\\
&= \chi(\lm - \xi(\blm))\chi(\mu - \xi(\bmu))\ \ph(x)\dotchi\ph(y)\\
& = \ph'(x) \dotchi \ph'(y).
\end{align*}
Further if $\gm\in \Gm$, $\lm\in \Lm$, $c\in C(\cB)^\gm$ and $x\in \cB^\lm$,
we have
\begin{align*}
\ph'(cx) &= \chi(\lm + \gm - \xi(\blm + \bar\gm))\ \ph(cx)\\
&= \chi(\lm + \gm - \xi(\blm ))\rho(c)\ \ph(x)\\
&= \chi(\gm)\chi(\lm  - \xi(\blm ))\rho(c)\ \ph(x)\\
&= \rho_\chi(c)\ph'(x).
\end{align*}
Thus $\ph'$ is a $\rho_\chi$-specialization and so $\cA_\chi$ is
a $\rho_\chi$-image of $\cA$.

(iii):  This follows from (ii) and the uniqueness of the $\rho_\chi$-image
(established in Corollary  \ref{cor:uniquecent}).
\end{proof}

\begin{proposition}
\label{prop:centim}
Suppose that
$\cB\in \classB$ and
$\cA$ is a central image of $\cB$. If $\cA'$ is a $\bLm$-graded
algebra, then $\cA'$ is a central image of $\cB$ if and only if
$\cA' \sim_\pi \cA$.
\end{proposition}

\begin{proof} We are given that $\cA$ is a $\rho$-image of $\cB$ for some
$\rho\in \algC$.

``$\Rightarrow$''  If $\cA'$ is a $\rho'$-image of
$\cB$  for some $\rho'\in \algC$,
then $\cA' \isombLm \cA_\chi$ for some $\chi\in \charG$ by
parts (i) and (iii) of Lemma \ref{lem:rhochi}. So $\cA\sim_\pi\cA'$.

``$\Leftarrow$" Suppose that $\cA^{\prime}\simeq\cA_{\chi}$ for
some $\chi\in \Hom(\Gamma,k^{\times})$.By Lemma \ref{lem:rhochi}(ii),
$\cA_{\chi}$, and hence $\cA^{\prime}$,  is a $\rho_\chi$-image of  $\cB$.
\end{proof}

\subsection{A characterization of central images}
\label{subsec:characterization}

The next
example shows that $\cA$ is a central image of $\loopA$
for $\cA\in \classA$.

\begin{example}
\label{ex:imloop}
Suppose that $\cA \in \classA$ and let $\cB = \loopA$.
Then
by \eqref{eq:CloopA}, we have
$C(\cB) = \spann_k\set{l_{1\ot z^\gm} \suchthat \gm\in \Gm}\isomLm k[\Gm]$,
and so we may define $\rho\in \algC$ by
\[\rho(l_{1\ot z^\gm})=1\]
for $\gm\in\Gm$.
We call $\rho$ the \textit{augmentation homomorphism}.
Define $\ph : \cB \to \cA$ by
\[\ph(u\ot z^\lm) = u\]
for $u\in \cA^\lm$ and $\lm\in \Lm$.
Then one checks easily that $\ph$ is a $\rho$-specialization.
Therefore, $\cA$ is a $\rho$-image of $\loopA$.
\end{example}

We can now characterize the central images
of a given $\cB\in \classB$ using loop algebras.

\begin{proposition}
\label{prop:characterize}
Suppose that $\cB\in\classB$ and
$\cA$ is a $\bLm$-graded algebra.
Then the following statements are equivalent:
\begin{itemize}
    \item[(a)] $\loopA \isomLm \cB$.
    \item[(b)] $\cA$ is a central image of $\cB$.
    \item[(c)] $\cA \isombLm \unim$ for some $\rho\in \algC$.
\end{itemize}
Moreover if (a), (b) or (c) hold, then
$\cA\in \classA$.
\end{proposition}

\begin{proof} If (a) holds then
$\cA \in \classA$  by  Proposition \ref{prop:loopGCS}.
Thus is suffices to show that (a), (b) and (c) are equivalent.
But (b) and (c) are equivalent
by the uniqueness of the $\rho$-image of $\cB$
for $\rho\in \algC$
(see Corollary \ref{cor:uniquecent}).  Moreover
``(b) $\Rightarrow$ (a)'' follows from
Proposition \ref{prop:spec}(iv).
Finally,
``(a) $\Rightarrow$  (b)''
follows from
Example~\ref{ex:imloop}.
\end{proof}

\begin{corollary}
\label{cor:simclosed}
Suppose
that $\cA\in \classA$ and $\cA'$ is a $\bLm$-graded algebra.
If $\cA\sim_\pi\cA'$ then $\cA'\in \classA$.
\end{corollary}

\begin{proof}  By Proposition \ref{prop:loopGCS},
$\loopA\in \classB$.  Then, by Proposition \ref{prop:characterize},
$\cA$ is a central image of $\loopA$.
Thus, since $\cA\sim_\pi\cA'$, we know by Proposition \ref{prop:centim}
that $\cA'$ is a central image of $\loopA$.  So, by Proposition \ref{prop:characterize},
$\cA'\in\classA$.
\end{proof}

\section{The correspondence}
\label{sec:correspondence}

Suppose again in this section that $\Gm$ is a subgroup of
an arbitrary abelian group
$\Lm$,
and that $\pi : \Lm \to \bLm$ is a  group epimorphism with  kernel $\Gm$.

\subsection{The correspondence theorem}
\label{subsec:main}

We can now combine the results from the previous sections to prove
our main theorem about the loop algebra construction.
This theorem  tells us that
the loop construction induces a correspondence
between  similarity classes of $\bLm$-graded algebras in $\classA$
and graded-isomorphism classes of $\Lm$-graded algebras in $\classB$.
The inverse correspondence
is induced by central specialization.

\begin{theorem}[Correspondence Theorem]
\label{thm:correspondence}
Let
$\Gm$ be a subgroup of $\Lm$ and let $\pi : \Lm \to \bLm$ be a group
epimorphism such that $\kerpi = \Gm$.
Let $\classA$ be the class of $\bLm$-graded algebras $\cA$
such that $\cA$ is central-simple as an algebra. Let
$\classB$ be the class of $\Lm$-graded algebras $\cB$ such that
$\cB$ is graded-central-simple, $C(\cB)$ is split and
$\Gm(\cB) = \Gm$.
For a
$\bLm$-graded  algebra $\cA$,
let
$\loopA = \sum_{\lm\in\Lm} \cA^\blm \ot z^\lm$, where $\blm = \pi(\lm)$ (see Definition~\ref{def:loop}).
\begin{itemize}
    \item[(i)] If $\cA\in \classA$, then $\loopA\in \classB$.
    \item[(ii)] If $\cB\in \classB$, then there exists $\cA\in \classA$ such
    that $\loopA \isomLm \cB$.  Moreover the $\bLm$-graded algebras $\cA$ with
    this property are precisely the central images of $\cB$.
    \item[(iii)] If $\cA,\cA'\in \classA$, then $\loopA \isomLm \loopAp$
    if and only if $\cA\sim_\pi \cA'$.
    \item[(iv)] If  $\cA\in \classA$, $\cB\in \classB$ and $\cB \isomLm
    \loopA$, then $\cA$ is finite dimensional if and only if $\cB$ is fgc
    (finitely generated as
    a module over its centroid).
\end{itemize}
\end{theorem}

\begin{proof}
Now (i) follows from Proposition \ref{prop:loopGCS}, and
(ii) follows from Proposition \ref{prop:characterize}.  To prove (iii),
suppose that $\cA,\cA'\in \classA$.  Then $\cA$ is a central image
of $\loopA$ by Proposition \ref{prop:characterize}.  So
\begin{align*}
\cA\sim_\pi\cA' &\iff \cA' \text{ is a central image of } \loopA
&(\text{by Proposition~\ref{prop:centim}})\\
&\iff \loopA \isomLm \loopAp &(\text{by Proposition~\ref{prop:characterize}}).
\end{align*}
Finally, (iv) follows from (ii) and  Proposition \ref{prop:spec}(ii).
\end{proof}

\begin{definition}
\label{def:correspond}
Let
$\Gm$, $\Lm$ and  $\pi : \Lm \to \bLm$ be as in the Correspondence Theorem.
If $\cA\in \classA$ and
$\cB \in \classB$, we say that $\cA$ and $\cB$ \textit{correspond}
(under $\pi$) if $\loopA \isomLm \cB$.
By  the theorem, $\cA$ and $\cB$
correspond under $\pi$ if and only if
$\cA$ is a central image of $\cB$ (relative to $\pi$).
\end{definition}

\begin{remark}
\label{rem:cat}
Theorem  \ref{thm:correspondence}
has a categorical formulation.  To describe this we define two categories
$\catA$ and $\catB$ and a functor $F: \catA \to \catB$. Although we
suppress this from the notation, $\catA$, $\catB$ and $F$ depend
on the group epimorphism $\pi: \Lm \to \bLm$ with kernel $\Gm$.
(The category $\catA$ and the functor $F$, but not the category $\catB$,
also depend on the choice of a fixed right inverse
$\xi$ for $\pi$ as in \S \ref{sec:central-spec}.)

The objects in $\catA$ are the $\bLm$-graded algebras in $\classA$. A
morphism in $\catA$ from $\cA$ to $\cA'$ is a pair $(\eta,\chi)$, where
$\eta$ is a $\bLm$-graded-isomorphism from $\cA$ to $\cA'_\chi$ and
$\chi\in \Hom(\Gm,k^\times)$.

The objects in $\catB$ are the pairs
$(\cB,\rho)$, where $\cB$ is a $\Lm$-graded algebra in $\classB$
and $\rho\in \Alg(C(\cB),k)$.
A morphism in $\catB$ from $(\cB,\rho)$ to $(\cB',\rho')$
is a $\Lm$-graded isomorphism $\psi : \cB \to \cB'$.

If $\cA$ is an object in $\catA$, we let
$F(\cA) = (\loopA,\rho)$, where
$\rho : C(\loopA)\to k$
is the augmentation homomorphism (see Example \ref{ex:imloop}).
If $(\eta,\chi)$
is a morphism from $\cA$ to $\cA'$ in $\catA$, we let
$F(\eta,\chi)$ be the map from $\loopA$ to $\loopAp$
given by
$x\textstyle \ot z^\lm \mapsto  \frac{1}{\rule{0em}{.7em}\chi(\lm-\xi(\blm))}\, \eta(x)\ot z^\lm$
for $x\in \cA^\lm$.

Then using the results (and their proofs) from
Sections \ref{subsec:propspec},
\ref{subsec:compare}
and \ref{subsec:characterization}
one can show that $F$ is an equivalence of categories.  We leave the details
of this to the  reader.

\begin{remark}
\label{rem:histcorr} Suppose  that $k$ is algebraically
closed of characteristic 0 and $\Lm$ is finite, and assume that the graded
algebras $\cA$ and $\cB$ are finite dimensional.   Then,
(i) and the  first statement of (ii) in Theorem \ref{thm:correspondence} were proved
by Bahturin, Sehgal and Zaicev in \cite[Theorem 7]{BSeZ}
(although the description of the loop algebra
used there is not the same as the one used  here).
\end{remark}

\end{remark}

\subsection{
A quantum torus example
}

\begin{example}
\label{ex:quantum}
Suppose that $m$ is a positive integer and that
$k$ contains a primitive $m^\text{th}$ root of unity $\zeta_m$.
Let $\Lm = \Ztwo$ and $\bLm = \bbZ/m\bbZ \oplus \bbZ/m\bbZ$
and let $\pi : \Lm \to \bLm$ be the natural map.

Let
$\boldq =\left[\begin{smallmatrix} 1&\zeta_m\\ \zeta_m^{-1} & 1 \end{smallmatrix}\right]$
and let $\cB = k_\boldq$ be the  \textit{quantum torus} determined by $\boldq$.
Thus, by definition, $\cB$ is the unital associative algebra determined
by the generators $x_1,x_1^{-1},x_2,x_2^{-1}$ subject to the relations
\[x_i x_i^{-1} =  x_i^{-1} x_i = 1 \text{  for  } i=1,2 \andd x_2x_1 = \zeta_m x_1x_2.\]
Then
$\cB = \bigoplus_{\lm\in \Lm} \cB^\lm$ is a
$\Lm$-graded algebra with $\cB^{(\ell_1,\ell_2)} = kx_1^{\ell_1}x_2^{\ell_2}$
for $(\ell_1,\ell_2)\in \Lm$.

We identify the centroid $C(\cB)$ of $\cB$ with the centre of $\cB$
(see Remark \ref{rem:centre})  in which case we have
\[C(\cB) = k[c_1^{\pm 1},c_2^{\pm 1}],\]
where $c_1 = x_1^m$ and $c_2 = x_2^m$.
Thus, $\cB$ is a graded-central-simple $\Lm$-graded algebra, $C(\cB)$ is split
and $\Gm(\cB) = m\bbZ\oplus m\bbZ$.  In other words,
$\cB\in {\mathfrak B(\Lm,m\bbZ\oplus m\bbZ)}$.

Next the elements of $\Alg(C(\cB),k)$ are the maps of the form
\[\rho_{a_1,a_2} : f(c_1,c_2)\mapsto f(a_1,a_2)\]
for $f(c_1,c_2)\in C(\cB)$, where $a_1,a_2\in k^\times$.
Let
\[\cA_{a_1,a_2} = \cB/\ker(\rho_{a_1,a_2})\cB\]
(as  in Example
\ref{ex:universal}) for $a_1,a_2\in k^\times$.  Then,
$\cA_{a_1,a_2}$ can be identified with the unital associative
 algebra determined by the generators
$y_1,y_2$ subject to the relations
\[y_1^m = a_1, \quad y_2^m = a_2 \andd y_2y_1 = \zeta_m y_1y_2.\]
That is, $\cA_{a_1,a_2}$ is the \textit{power norm residue algebra}
over $k$ determined by $a_1$, $a_2$, $m$ and $\zeta_m$ (see \cite[\S 11]{D}).
Moreover, the $\bLm$-grading on $\cA_{a_1,a_2}$
is determined by the conditions
\[\deg(y_1) = (\bar 1,\bar 0) \andd \deg(y_2) = (\bar 0,\bar 1) .\]
Thus, by Proposition \ref{prop:characterize},
the central images of $\cB$ are precisely the power norm residue algebras
$\cA_{a_1,a_2}$, $a_1,a_2\in k^\times$ (with the indicated $\bLm$-gradings).
Consequently, these algebras are precisely the $\bLm$-graded algebras
in $\classA$ that correspond to $\cB$ under $\pi$.

It follows from Theorem \ref{thm:correspondence}(ii), that
we have
\[L_\pi(\cA_{a_1,a_2}) \isomLm \cB\]
for $a_1,a_2\in k^\times$.  Consequently
(see \eqref{eq:mloopisaloop}), we have
\begin{equation}
\label{eq:misomquant}
M_{(m,m)}(\cA_{a_1,a_2},\sg_1,\sg_2) \isomLm\cB,
\end{equation}
where $\sg_1$, $\sg_2$ are the automorphisms of $\cA_{a_1,a_2}$
determined by the conditions:
\[\sg_1(y_1) = \zeta_m y_1,\quad \sg_1(y_2) = y_2,\quad \sg_2(y_1) = y_1 \andd \sg_2(y_2) = \zeta_m y_2.\]
\end{example}

\medskip
\begin{remark}
The  fact that the quantum torus $k_\boldq$ in
Example \ref{ex:quantum} is a multiloop algebra was  previously observed
in \cite[Example 9.8]{ABP}, where it was shown directly
that $\cB$ is a multiloop algebra
based on $\cA_{1,1}$ (the algebra of $m\times m$-matrices over $k$).  The isomorphism
\eqref{eq:misomquant} shows more generally that $\cB$ is a multiloop algebra
based on any power norm residue algebra $\cA_{a_1,a_2}$ and it places
this fact in the much more general context of the correspondence theorem.
\end{remark}

\section{Multiloop realization of graded algebras}
\label{sec:mult}

In this section  we prove our main results about
multiloop realizations of graded-central-simple algebras
(see Theorem \ref{thm:realmult} and Corollary
\ref{cor:multreal}).
Throughout the section we assume that
$n$ is an integer $\ge 1$, $\Lm$ is a free abelian group of rank $n$ and
$k$ is an algebraically closed field
of characteristic zero.

In view of our assumptions on $k$,
$k$ contains a primitive root $\ell^\text{th}$
of unity $\zeta_\ell$ for all positive integers $\ell$. We assume that
we have made a fixed compatible choice of these roots of unity
in the sense that
\begin{equation}
\label{eq:compat1}
\zeta_{m \ell}^m =
\zeta_\ell
\end{equation}
for all $\ell, m\ge 1$.  (This is always possible.)
We use these roots of unity in the construction
of multiloop algebras.

\subsection{$(\boldmp,\boldm)$-admissible matrices}
We  begin by  describing some terminology that will be useful
in the study of multiloop algebras.

\begin{definition}
\label{def:related}
Let
$\boldm = (m_1,\dots,m_n)$ and $\boldmp
= (m'_1,\dots,m'_n)$ be $n$-tuples of positive integers.
Let $D_\boldm = \diag(m_1,\dots,m_n)$ and $D_\boldmp = \diag(m'_1,\dots,m'_n)$.

(i) If $P$ is an $n\times n$-matrix with rational entries, we
define the $(\boldmp,\boldm)$-\textit{transpose} of $P$ to be the
$n\times n$-matrix
\[
D_\boldmp P^t D_\boldm^{-1} = (\frac{m'_i}{m_j} p_{ji}),
\]
where $P^t$ denotes the (usual) transpose of $P$.

(ii) Recall that $\GL_n(\bbZ)$ is the group of all $n\times n$
matrices with integer entries and determinant $\pm 1$.
If $P \in\GL_n(\bbZ)$, we say that $P$ is $(\boldmp,\boldm)$-\textit{admissible}
if the $(\boldmp,\boldm)$-transpose of $P$ is in $\GL_n(\bbZ)$.
Note that
if $P$ is $(\boldmp,\boldm)$-admissible with $(\boldmp,\boldm)$-transpose $Q$,
then $Q$ is $(\boldm,\boldmp)$-admissible with $(\boldm,\boldmp)$-transpose $P$
and $P^{-1}$ is $(\boldm,\boldmp)$-admissible with $(\boldm,\boldmp)$-transpose~$Q^{-1}$.
\end{definition}

The following lemma gives a useful characterization
of  $(\boldmp,\boldm)$-admissible matrices.  In  this lemma
(and later) we identify $P\in \GL_n(\bbZ)$
with the automorphism of $\Zn$ given by
\begin{equation}
\label{eq:identGL}
\boldl \overset{P}{\longrightarrow} \boldl P^t.
\end{equation}
(We use the right transpose action since we are viewing elements of $\Zn$
as row vectors.)

\begin{lemma}
\label{lem:related}  Let $\boldm=(m_{1},\ldots,m_{n})$ and
$\boldmp=(m_{1}^{\prime},\ldots,m_{n}^{\prime})$ be
$n$-tuples of positive integers and let $P\in GL_{n}(\bbZ)$.
Set $\bLm=\bbZ/(m_{1})\oplus\ldots\oplus\bbZ/(m_{n})$ and
$\bLm^{\prime}=\bbZ/(m_{1}^{\prime})\oplus\ldots\oplus\bbZ
/(m_{n}^{\prime})$. Then $P$ induces an isomorphism of
$\bLm^{\prime}$ onto $\bLm$  if and only if $P$
is $(\boldmp,\boldm)$-admissible.
\end{lemma}

\begin{proof}
We note that the kernel of natural homomorphism
$\Zn\rightarrow \bLm$ is
$\Zn D_{\boldm}$. Thus, $P$ induces an isomorphism of
$\bLm^{\prime}$ onto $\bLm$ if and only if
$\Zn D_{\boldm^{\prime}}P^{t}=\bbZ^{n}D_{\boldm}$,
or
equivalently $\bbZ^{n}D_{\boldm^{\prime}}P^{t}D_{\boldm}
^{-1}=\bbZ^{n}$. But the last condition is equivalent to
$D_{\boldm^{\prime}}P^{t}D_{\boldm}^{-1}\in GL_{n}(\bbZ)$.
\end{proof}

\subsection{Properties of multiloop algebras}
\label{subsec:propmulti}

In this subsection, we prove two basic propositions
about multiloop  algebras.

\begin{proposition}
\label{prop:mult1}
 Let $\boldm=(m_{1},\ldots,m_{n})$
and $\boldm^{\prime}=(m_{1}^{\prime},\ldots,m_{n}^{\prime})$ be
$n$-tuples of positive integers and suppose that $P=(p_{ij})\in
GL_{n}(\bbZ)$ is
$(\boldm^{\prime},\boldm)$-admissible with
\begin{equation}
Q:=D_{\boldm^{\prime}}P^{t}D_{\boldm}^{-1}\in GL_{n}(\bbZ
).\label{eq:Q}
\end{equation}
Let $\sg_{1},\ldots,\sg_{n}$ be commuting automorphisms of an
algebra
$\cA$ such that $\sg_{i}^{m_{i}}=1$ for $1\leq i\leq n$, and let
\[
\sg_{i}^{\prime}=
{\textstyle\prod\nolimits_{j=1}^{n}}
\sg_{j}^{p_{ji}}
\]
for $1\leq i\leq n$. Then
$\sg_{1}^{\prime},\ldots,\sg_{n}^{\prime}$ are commuting
automorphisms of $\cA$ such that $\sg_{i}^{\prime
m_{i}^{\prime}}=1$ for $1\leq i\leq n$. Moreover, we have
\begin{equation}
M_{\boldm^{\prime}}(\cA,\sg_{1}^{\prime},\ldots,\sg
_{n}^{\prime})\simeq_{\bbZ^{n}}M_{\boldm}(\cA,\sg
_{1},\ldots,\sg_{n})_{R}\label{eq:graded isom}
\end{equation}
where $R=Q^{-1}\in GL_{n}(\bbZ)$, and hence
\begin{equation}
M_{\boldm^{\prime}}(\cA,\sg_{1}^{\prime},\ldots,\sg
_{n}^{\prime})\simeq_{ig}M_{\boldm}(\cA,\sg_{1},\ldots
,\sg_{n}).\label{eq:isograded}
\end{equation}
(Here we are using the notation in Definitions
\ref{def:isomorphism}(ii) and \ref{def:regrade},
as well as the identification $\Aut(\Zn) = \GL_n(\bbZ)$
in \eqref{eq:identGL}.)
\end{proposition}

\begin{proof}
As in Remark \ref{rem:mloop}(iii), set
$G=\langle \sg_{1},\ldots,\sg_{n}\rangle$
and
$\sg^{\boldl}=
{\textstyle\prod\nolimits_{i=1}^{n}}
\sg_{i}^{l_{i}}\in G$
for $\boldl=(l_{1},\ldots l_{n})\in\bbZ^{n}$. Similarly, set
\begin{equation}
\sg^{\prime\boldl}=
{\textstyle\prod\nolimits_{i=1}^{n}}
\sg_{i}^{\prime l_{i}}=
{\textstyle\prod\nolimits_{i=1}^{n}}
{\textstyle\prod\nolimits_{j=1}^{n}}
\sg_{j}^{l_{i}p_{ji}}=
{\textstyle\prod\nolimits_{j=1}^{n}}
\sg_{j}^{
{\textstyle\sum\nolimits_{i=1}^{n}}
l_{i}p_{ji}}=\sg^{\boldl P^{t}}.\label{eq:powers}
\end{equation}
Since $P$ is invertible it follows that
we also have
$G=\langle \sg'_{1},\ldots,\sg'_{n}\rangle$.

Now our assumption that $\sg_{i}^{m_{i}}=1$ for $1\leq i\leq n$ is equivalent to
$\sg ^{\bbZ^{n}D_{\boldm}}=1$. But we have
\[
\sg^{\prime\bbZ^{n}D_{\boldm^{\prime}}}=\sg^{\bbZ
^{n}D_{\boldm^{\prime}}P^{t}}=\sg^{\bbZ^{n}QD_{\boldm}}=1,
\]
and so $\sg_{i}^{\prime m_{i}^{\prime}}=1$ for $1\leq i\leq n$.

It remains to prove (\ref{eq:graded isom}) (since (\ref{eq:isograded}) follows
from (\ref{eq:graded isom})). For this purpose, we let $\cB
=M_{\boldm}(\cA,\sg_{1},\ldots,\sg_{n})$ and $\cB
^{\prime}=M_{\boldm^{\prime}}(\cA,\sg_{1}^{\prime}
,\ldots,\sg_{n}^{\prime})$, in which case we must prove that $\cB
_{R}\simeq_{\bbZ^{n}}\cB^{\prime}$ or equivalently
\begin{equation}
\cB\simeq_{\bbZ^{n}}\cB_{Q}^{\prime}\label{eq:B version}
\end{equation}
(since $Q=R^{-1}$).

Let $S=k[z_{1}^{\pm1},\ldots,z_{n}^{\pm1}]$. Then
$\mathcal{A\otimes}_{k}S$ is a $\bbZ^{n}$-graded algebra
(with the grading determined by the natural grading on $S$), and
both $\cB$ and $\cB^{\prime}$ are graded
subalgebras of $\mathcal{A\otimes}_{k}S$. Next we define $\gamma
_{Q}\in \Aut_{k}(S)$ by $\gamma_{Q}(z^{\boldl})=z^{\boldl\, Q^{t}}$
where
$z^{\boldl}=z_{1}^{l_{1}}\ldots z_{n}^{l_{n}}$ for $\boldl
=(l_{1},\ldots l_{n})$. Then $1\otimes\gamma_{Q}$ is a $\bbZ^{n}
$-graded algebra isomorphism of $\mathcal{A\otimes}_{k}S$ onto
$(\mathcal{A\otimes}_{k}S)_{Q}$. Consequently, to prove
\eqref{eq:B version}, it suffices to show that $1\otimes\gamma_{Q}$ maps $\cB$ onto
$\cB^{\prime}$.

Now let
$\bar{\Lm}=\bbZ/(m_{1})\oplus\ldots\oplus\bbZ/(m_{n})$
and $\bar{\Lm}^{\prime}=\bbZ/(m_{1}^{\prime})\oplus\ldots
\oplus\bbZ/(m_{n}^{\prime})$. We regard $\cA$ as a
$\bar{\Lm}$-graded algebra with the grading determined by $\sg
_{1},\ldots,\sg_{n}$ (see Definition 3.2.1). Let $\cA^{\prime}$
denote the algebra $\cA$ with the $\bar{\Lm}^{\prime}$-grading
determined by $\sg_{1}^{\prime},\ldots,\sg_{n}^{\prime}$. Then
$\cB=
{\textstyle\sum\nolimits_{\boldl\in\bbZ^{n}}}
\cA^{\boldlb}\otimes z^{\boldl}$ and so
\[
(1\otimes\gamma_{Q})\cB=
{\textstyle\sum\nolimits_{\boldl\in\bbZ^{n}}}
\cA^{\boldlb}\otimes z^{\boldl\, Q^{t}}.
\]
On the other hand, $\cB^{\prime}=
{\textstyle\sum\nolimits_{\boldl\in\bbZ^{n}}}
\cA^{\prime\boldlb}\otimes z^{\boldl}$, and so it suffices to
show that
\begin{equation}
\label{eq:l and lQ}
\cA^{\boldlb}=\cA^{\prime\overline{\boldl\, Q^{t}}}
\end{equation}
for  $\blm\in\bLm$.
Since $\zeta_{ml}^{m} =\zeta_{l}$ for all $m,l\geq1$,
$\zeta(\frac{m}{n}) := \zeta_{n}^{m}$
defines a homomorphism $(\mathbb{Q},+)\rightarrow k^{\times}
$. Now
\[
\cA^{\boldlb}=\{x\in\cA\mid gx=\chi_{\boldlb}(g)x\text{ for all }g\in G\},
\]
where $\chi_{\boldlb}$ is the character on $G$ with $\chi
_{\boldlb}(\sg_{i})=\zeta_{m_{i}}^{l_{i}}=\zeta(\frac{l_{i}}
{m_{i}})$. One checks  that $\chi_{\boldlb}(\sg^{\mathbf{k}}
)=\zeta((\boldl D_{\boldm}^{-1})\cdot\mathbf{k})$ where $\cdot$ is
the
usual dot product. Thus to prove \eqref{eq:l and lQ} it is enough to show that
\[
\chi_{\boldlb}=\chi_{\overline{\boldl\, Q^{t}}}^{\prime}
\]
where $\chi_{\overline{\boldl\, Q^{t}}}^{\prime}(\sg^{\prime\mathbf{k}
})=\zeta((\boldl Q^{t}D_{\boldm^{\prime}}^{-1})\cdot\mathbf{k})$.
Since $\sg^{\prime\mathbf{k}}=\sg^{\mathbf{k}P^{t}}$ by
(\ref{eq:powers}), the result follows from
\[
(\boldl D_{\boldm}^{-1})\cdot (\mathbf{k}P^{t})=(\boldl D_{\boldm
}^{-1}P)\cdot\mathbf{k}=(\boldl Q^{t}D_{\boldm^{\prime}}^{-1}
)\cdot\mathbf{k}.\qedhere
\]
\end{proof}

We  next use the Correspondence Theorem to prove a second proposition about multiloop algebras.

\begin{proposition}
\label{prop:mult2}\
\begin{itemize}
    \item[(i)] Suppose that
$\boldm= (m_1,\dots,m_n)$ is an $n$-tuple of positive integers,
$\cA$ is a central-simple (ungraded) algebra
and $\sg_1,\dots,\sg_n$ are commuting algebra automorphisms
of $\cA$ such that $\sg^{m_i} = 1$ for each $i$.  Then
$\mloopA$ is a graded-central-simple $\Zn$-graded algebra
whose central grading group is given by
\begin{equation}
    \label{eq:Gmmult}
    \Gm(\mloopA) = m_1\bbZ\oplus \dots \oplus m_n\bbZ.
\end{equation}
\item[(ii)] Suppose that  $\boldm$, $\cA$ and $\sg_1,\dots,\sg_n$ are as in (i)
and $\boldm'$, $\cA'$ and $\sg'_1,\dots,\sg'_n$ are as in (i).  Then
$\mloopA \simeq_{\Zn} \mloopAp$ if and only if $\boldm = \boldmp$
and there exists an algebra isomorphism
$\ph : \cA \to \cA'$ such that
\begin{equation}
\label{eq:phconj}
\ph \sg_j \ph^{-1} = \sg'_j
\end{equation}
for $1\le j\le n$.
\end{itemize}
\end{proposition}

\begin{proof}
(i) Let $\cB = \mloopA$ and $\Lm = \Zn$.  Then, as we saw in
Definition \ref{def:multiloop}, we have $\cB = \loopASig$,
where $\bLm = \bbZ/(m_1)\oplus \dots \oplus
\bbZ/(m_n)$, $\pi : \Lm \to \bLm$ is the natural map defined by
\eqref{eq:natpi}, and
$\Sigma$ is the $\bLm$-grading on $\cA$ determined by the
automorphisms $\sg_1,\dots,\sg_n$. By part (i) of the Correspondence Theorem,
$\cB$ is graded-central-simple with
central grading group $\kerpi$. Since $\kerpi = m_1\bbZ\oplus \dots
\oplus m_n\bbZ$, we have   \eqref{eq:Gmmult}.

(ii): Let $\cB = \mloopA$,  $\cB' = \mloopAp$ and $\Lm = \Zn$.
If $\cB \isomLm \cB'$, then $\Gm(\cB) = \Gm(\cB')$ and so we have $\boldm = \boldmp$
by \eqref{eq:Gmmult}.
Consequently, for the rest of the proof of both directions in (ii), we can and do
assume that  $\boldm = \boldmp$.

Let $\bLm$ and $\pi$ be as in the proof of (i), and let $\Sigma$
(resp.~$\Sigma'$) be the $\bLm$-grading on $\cA$ (resp~$\cA'$) determined by the
automorphisms $\sg_1,\dots,\sg_n$ (resp.~$\sg'_1,\dots,\sg'_n$).
Then
$\cB = \loopASig$ and $\cB' = L_\pi(\cA',\Sigma')$.
Hence, by part (iii) of the Correspondence Theorem, it follows that
$\cB \isomLm \cB'$ if and only if $(\cA,\Sigma) \sim_\pi
(\cA',\Sigma')$.  Furthermore, since $k$  is algebraically closed,
we have by
Remark \ref{rem:twist}(iii)
that $(\cA,\Sigma)
\sim_\pi (\cA',\Sigma')$ if and only if $(\cA,\Sigma)\isombLm
(\cA',\Sigma')$.   Consequently, it suffices to show that the
$\bLm$-graded algebra isomorphisms from $(\cA,\Sigma)$ to
$(\cA',\Sigma')$ are precisely the algebra isomorphisms $\ph : \cA
\to \cA'$ that satisfy \eqref{eq:phconj}.  This is easily  checked.
\end{proof}

\subsection{Multiloop realizations}
\label{subsec:realmult}

\begin{definition}
\label{def:classBfi}
Let $\classBfi$ be the class of $\Lm$-graded
algebras $\cB$ such that $\cB$ is graded-central-simple and
$\Gm(\cB)$ has finite index in $\Lm$.
Note that
since
$\Lm$ is free of finite rank
(by assumption),
it follows from Lemma \ref{lem:centsplit2} that
any graded algebra $\cB$
in $\classBfi$ has split centroid.
Hence
\[\classBfi = \cup_\Gm  \classB,\]
where the class union $\cup_\Gm$ runs over all subgroups
$\Gm$ of finite index in~$\Lm$.
\end{definition}

Our next main result gives multiloop realizations and isomorphism conditions
for all graded algebras in $\classBfi$.
In this theorem we use the notion of \textit{isograded-isomorphism} and the notation $\ig$
described
in Definition \ref{def:isomorphism}(ii).

\begin{theorem}[Realization Theorem]
\label{thm:realmult}
Suppose that $k$ is an algebraically closed field of characteristic 0.
\begin{itemize}
\item[(i)]  Suppose that $\cB$ is a $\Lm$-graded algebra,
where $\Lm$ is a free abelian group of rank $n\ge 1$.
 Then
    $\cB\in \classBfi$ if and only if $\cB$ is isograded-isomorphic
    to $\mloopA$ for some central-simple (ungraded) algebra $\cA$,
    some $n$-tuple of positive integers $\boldm= (m_1,\dots,m_n)$
    and some  sequence $\sg_1,\dots,\sg_n$ of commuting finite order algebra automorphisms
    of $\cA$ such that $\sg^{m_i} = 1$ for all $i$.

\item[(ii)]
Let $\Lm = \Zn$.
Suppose $\cB= \mloopA$, where $\cA$, $\boldm$, and $\sg_1,\dots,\sg_n$ are as in
(i), and
suppose
$\cB' = \mloopAp$, where $\cA'$, $\boldm'$, and
$\sg'_1,\dots,\sg'_n$ are as in (i). Then $\cB\ig \cB'$ if and only
if there exists a matrix $P = (p_{ij}) \in\GL_n(\bbZ)$ and an algebra
isomorphism $\ph : \cA\to \cA'$ such that $P$ is
$(\boldm',\boldm)$-admissible and
\begin{equation}
\label{eq:conjcond}
      \sg'_j  = \ph \left(\textstyle\prod_{i=1}^n \sg_i^{p_{ij}}\right) \ph^{-1}
      \end{equation}
for $1\le j\le n$.
Moreover, if
$\langle \supp_\Zn(\cB)\rangle = \Zn$ and
$\langle \supp_\Zn(\cB')\rangle = \Zn$, then
$\cB\ig \cB'$ if and only
if there exists a matrix $P = (p_{ij}) \in\GL_n(\bbZ)$ and an algebra
isomorphism $\ph : \cA\to \cA'$ such that \eqref{eq:conjcond}
holds for $1\le j\le n$.

\item[(iii)]
Suppose $\cB\in
\classBfi$,
where $\Lm$ is a free abelian group of rank $n\ge 1$,
and suppose
$\cB \ig \mloopA$,
where $\cA$, $\boldm$, and
$\sg_1,\dots,\sg_n$ are as in~(i).
Then $\cA$ is finite
dimensional if and only if $\cB$ is fgc.
\end{itemize}
\end{theorem}

\begin{proof}
(i):  The implication ``$\Leftarrow$'' follows from
Proposition \ref{prop:mult2}(i).  To prove ``$\Rightarrow$'', let
$\cB \in \classBfi$.  Then $\cB \in \classB$ for some subgroup
$\Gm$ of $\Lm$ of finite index in~$\Lm$.
Furthermore, by the fundamental
theorem of finitely generated abelian groups, there
exists a $\bbZ$-basis  $\set{\lm_1,\dots,\lm_n}$ of $\Lm$
and an $n$-tuple $\boldm = (m_1,\dots,m_n)$ of positive integers such that
$m_{i} | m_{i+1}$ for $1\le i \le n-1$
and
$\Gm = \langle m_1\lm_1,  \dots , m_n\lm_n\rangle$.
We then identify $\Lm = \Zn$ is such a way that
$\set{\lm_1,\dots,\lm_n}$ is the standard basis.  (We can do this since
we are working up to isomorphism of the grading groups.)
Let
$\bLm = \bbZ/(m_1)\oplus \dots \oplus \bbZ/(m_n)$, and let $\pi : \Lm \to \bLm$ be
the natural map.  Then $\pi$ is an epimorphism with kernel $\Gm$
and so by Theorem \ref{thm:correspondence}(ii) there exists $(\cA,\Sigma)\in \classA$
such that $\cB \isomLm \loopASig$.  (Here as in Definition \ref{def:multiloop}, it is convenient
to not abbreviate the graded algebra $(\cA,\Sigma)$ as $\cA$.)
Now since $\Sigma = \set{\cA^\blm}_{\blm\in\bLm}$
is a $\bLm$-grading there exist unique algebra automorphisms
$\sg_1,\dots,\sg_n$ of $\cA$ such that
\[
\cA^{(\modell_1,\dots,\modell_n)} = \set{u\in \cA \suchthat
\sg_j u = \zeta_{m_j}^{\ell_j} u \text{ for } 1\le j\le n}
\]
for $(\ell_1,\dots,\ell_n)\in \Zn$.  Then
$\sg_1,\dots,\sg_n$ is a sequence of commuting algebra automorphisms of $\cA$
such that $\sg_i^{m_i} = 1$ for all $i$.  Furthermore,
as we saw in Definition \ref{def:multiloop}, we have $\loopASig = \mloopA$.
Thus $\cB \isomLm \mloopA$.

(ii): By Remark \ref{rem:regrade}, we have
\[\cB \ig \cB' \iff \cB' \simeq_\Zn \cB_\nu\]
for some $\nu\in \Aut(\Zn)$.

Suppose for the moment that $\cB' \simeq_\Zn \cB_\nu$, where $\nu\in
\Aut(\Zn)$. Then $\Gm(\cB') = \nu^{-1}(\Gm(\cB))$ (by Remark
\ref{rem:centgroup}) and so $\nu(\Gm(\cB')) = \Gm(\cB)$.
So if we identify $\nu$ with  a matrix  $R\in \GL_n(\bbZ)$
(as in \eqref{eq:identGL}), we  have
$(\Zn D_\boldmp) R^t = \Zn D_\boldm$ by
\eqref{eq:Gmmult}.  Therefore
$R$ is $(\boldm',\boldm)$-admissible by
Lemma  \ref{lem:related}.

Consequently
\[\cB \ig \cB' \iff \cB' \simeq_\Zn \cB_R\]
for some $(\boldm',\boldm)$-admissible matrix $R\in \GL_n(\bbZ)$.  The first
statement in (ii) now follows from Proposition
\ref{prop:mult1} and  Proposition \ref{prop:mult2}(ii).

To prove the second statement in (ii),  suppose that
$\langle \supp_\Zn(\cB)\rangle = \Zn$ and
$\langle \supp_\Zn(\cB')\rangle = \Zn$.
It suffices to show that if $\ph : \cA \to \cA'$ is an algebra isomorphism
and  $P = (p_{ij}) \in\GL_n(\bbZ)$ satisfies \eqref{eq:conjcond} for $1\le j\le n$,
then $P$ is necessarily $(\boldmp,\boldm)$-admissible.

Now \eqref{eq:conjcond}  implies that
\[
\sg^{\prime\boldl}=\varphi\sg^{\boldl P^{t}}\varphi^{-1}
\]
for $\boldl\in\bbZ^{n}$. Thus,
there  is an isomorphism of
$\langle \sg'_{1},\ldots,\sg'_{n}\rangle$ onto
$\langle \sg_{1},\ldots,\sg_{n}\rangle$ such that
$\sg^{\prime\boldl }\rightarrow\sg^{\boldl P^{t}}$.
Hence, by Lemma \ref{lem:amplesupport}, there is an isomorphism from
$\bLm'$ onto $\bLm$ such that $\boldlb \mapsto \boldlb P^t$.
So, by Lemma \ref{lem:related},
$P$ is $(\boldm^{\prime},\boldm)$-admissible.

(iii):  This  follows from  Theorem \ref{thm:correspondence}(iv).
\end{proof}

\begin{definition}  If $\cB$ is a $\Lm$-graded algebra, $\cA$ is an algebra
and $\boldm$ is an $n$-tuple of positive integers, we say that $\cB$ has a \textit{multiloop
realization} based on $\cA$ and relative to $\boldm$ if there exist
commuting automorphisms $\sg_1,\dots,\sg_n$ of $\cA$ such that
$\cB \ig \mloopA$.
\end{definition}

\begin{remark}
Suppose that $\cB\in \classBfi$. Then by the Realization
Theorem (i), $\cB$ has a multiloop algebra realization
based on some central-simple algebra $\cA$ and relative to some
$\boldm$.   We note  now that more can be said about the choice of $\boldm$ and~$\cA$.

(i)  By the proof of part (i) of the Theorem, we may choose $\boldm$
with the additional property that
$m_{i} | m_{i+1}$ for $1\le i \le n-1$.
In that case, $\boldm$ is uniquely determined by $\cB$.
Indeed if $\cB$ has a multiloop algebra realizations relative to both
$\boldm$ and $\boldmp$,
then by part (ii) of the Theorem
there exists an $(\boldm',\boldm)$-admissible matrix.
But then by Lemma \ref{lem:related},
$\bbZ/(m_1)\oplus\dots\oplus\bbZ/(m_n) \simeq \bbZ/(m'_1)\oplus\dots\oplus\bbZ/(m'_n)$.
Hence if $m_{i} | m_{i+1}$ and $m'_{i} | m'_{i+1}$ for all $i$, we have
$\boldm = \boldmp$ by the fundamental theorem for finitely generated abelian groups.

(ii)  By part (ii) of the theorem, the ungraded algebra
 $\cA$ is uniquely determined up to isomorphism by $\cB$.
Moreover,  by the proof of part (i) of the theorem and by
Theorem \ref{thm:correspondence}(ii),
we may take $\cA$  to be any central image of $\cB$ (forgetting
the grading on $\cA$).
\end{remark}

The Realization Theorem has the following corollary
about multiloop realizations based on finite dimensional  simple algebras.

\begin{corollary}
\label{cor:multreal}
Let $\cB$ be a $\Lm$-graded algebra.
Then $\cB$ has a multiloop realization based on a finite dimensional  simple
algebra if and only if   $\cB$ is graded-central-simple, $\cB$ is fgc,
and $\Gm(\cB)$ has finite index in~$\Lm$.  Moreover,
if there exists a positive integer $\ell$ such that
$\ell\Lm \subseteq \supp_\Lm(\cB)$, then
$\cB$ has a multiloop realization based on a finite dimensional simple
algebra if and only if  $\cB$ is graded-central-simple and fgc.
\end{corollary}

\begin{proof}  Since  any finite dimensional simple algebra
over the algebraically closed field $k$ is central simple,
the first statement follows from
from parts (i) and (iii) of the  Realization theorem.
The second statement then follows from
the implication ``(a) $\Rightarrow$ (d)'' in Proposition \ref{prop:bornmodfree}.
\end{proof}

\begin{remark}  In  view of Proposition \ref{prop:bornmodfree},
the first statement in the preceding corollary can
be stated  alternatively as follows:
$\cB$ has a multiloop realization based on a finite dimensional  simple
algebra if and only if
$\cB$ is graded-central-simple, $\dim \cB^\lm < \infty$ for all $\lm\in \Lm$,
and $\Gm(\cB)$ has finite index  in~$\Lm$.
\end{remark}

\section{Some classes of tori}
\label{sec:tori}

The quantum torus discussed in Example \ref{ex:quantum}
is an example of what is called an associative torus.  (In fact this example
explains the use of the term associative torus.)  Furthermore
there are nonassociative analogs of these algebras called
alternative tori and Jordan tori, and there are Lie algebra
analogs called Lie tori.  When $\Lm$ is free of finite rank and $\charr(k) = 0$,
centreless Lie tori (Lie tori with trivial centre)
play a basic role in the theory of extended affine Lie algebras
because they appear as centreless cores of EALA's \cite{Y2,N}.  Furthermore,
associative, alternative and Jordan tori are also of great importance
in this context because they arise as coordinate algebras
of Lie tori of type $A_\ell$ (see \cite{BGK},
\cite{BGKN} and \cite{Y1}).  Thus an understanding
of these classes of tori is very important in the theory of EALA's.
In this section, we consider associative, alternative and Jordan tori,
leaving Lie tori for a  separate paper.

Initially,
we suppose only that
$k$ is a field and $\Lm$ is an abelian group.

\subsection{Definitions}
\label{subsec:torusdef}
\begin{definition}  Suppose that
$\cB$ is a unital associative, alternative or Jordan $\Lm$-graded algebra.
(We assume that $k$ has characteristic $\ne 2$ in the Jordan case.)
Then $\cB$ is said to be an \textit{associative, alternative or Jordan
$\Lm$-torus} respectively if
\begin{itemize}
\item[(a)] For each $\lm\in \supp_\Lm(\cB)$, $\cB^\lm$ is spanned
by an invertible element of $\cB$.
\item[(b)]  $\langle \supp_\Lm(\cB) \rangle = \Lm$.
\end{itemize}
(The interested reader can consult \cite[\S 10.3 and \S 14.2]{ZSSS}
for the basic
facts about invertibility in alternative and Jordan algebras.)
\end{definition}

\begin{proposition}
\label{prop:toriprop}
Suppose that $\cB$ is an associative, alternative or Jordan $\Lm$-torus.
Then $\cB$ is a graded-central-simple algebra and there exists
a positive integer $\ell$ such that
$\ell\Lm \subseteq \supp_\Lm(\cB)$.  In fact,
we may take $\ell= 1$ for associative and alternative tori,
and $\ell = 2$ for Jordan tori.
\end{proposition}

\begin{proof}  It is clear from the definition
that $\cB$ is graded-simple.  Also
$\dim_k(\cB_\lm) = 1$ for $\lm\in \supp_\Lm(\cB)$,
and so by Lemma \ref{lem:automaticcent}, $\cB$ is graded-central-simple.

It remains to check the last statement in the proposition.  Let
$T = \supp_\lm(\cB)$.
If $\cB$ is an associative or alternative torus,
then $T$ is a subgroup of $\Lm$ (since
the product of two invertible elements is invertible)
and so $T = \Lm$.  Suppose next that $\cB$ is a Jordan torus.
Then, $0\in T$, $-T = T$ and $T+2T \subseteq T$ (see
\cite[Lemma 3.5]{Y1}).  Since $\Lm = \langle T \rangle$,
it follows that
$T + 2\Lm \subseteq T$ and so
$2\Lm\subseteq T$ as desired.
\end{proof}

\subsection{Multiloop realization of tori}

For the rest of the paper,
we assume again that $k$ is algebraically closed of characteristic
$0$ and $\Lm$ is free abelian of finite rank~$\ge 1$.

We have the following application of our results:

\begin{theorem}
\label{thm:realtorus}  Suppose that $k$ is an algebraically closed
field of characteristic~0, and $\Lm$ is a free abelian group of finite rank  $\ge 1$.
Suppose that $\cB$ is an associative,
alternative or Jordan $\Lm$-torus.  Then
$\cB$ has a multiloop realization based on a finite dimensional
 simple associative, alternative or Jordan algebra $\cA$ if and only if
$\cB$ is~fgc.
\end{theorem}

\begin{proof}  This follows from Corollary \ref{cor:multreal} and
Proposition \ref{prop:toriprop}.  (See also Remark \ref{rem:loopgen}(iii).)
\end{proof}

\begin{remark}
\label{rem:realtorus}
Let
$\Lm = \Zn$
and
$\cB = \mloopA$, where
$\boldm = (m_1,\dots,m_n)$ is a sequence of positive integers,
$\cA$ is a finite dimensional  simple associative, alternative or Jordan
algebra, and  $\sg_1,\dots,\sg_n$ is a sequence
of commuting automorphisms of $\cA$ with $\sg_i^{m_i} = 1$ for all $i$.
Then, $\cB$ is an associative, alternative
or Jordan $\Lm$-torus respectively if and only if
the following conditions hold:
\begin{itemize}
\item[(a)]  The simultaneous eigenspaces  in $\cA$ for
the automorphisms $\sg_1,\dots,\sg_n$ are each spanned by
an invertible element of $\cA$.
\item[(b)]
$\order{\langle \sg_1,\dots,\sg_n\rangle}  = m_1\cdots m_n$.
\end{itemize}
We omit the verification of this observation which is
straightforward using Lemma \ref{lem:amplesupport}.
\end{remark}

In view of Theorem \ref{thm:realtorus} and Remark \ref{rem:realtorus}
(as well as the isomorphism condition
in the last sentence of part (ii) of the Realization Theorem),
the study of fgc associative, alternative or Jordan $\Lm$-tori
is equivalent to the study of sequences of commuting automorphisms of finite dimensional
algebras satisfying the conditions (a) and (b) in Remark \ref{rem:realtorus}.

\begin{remark}  Proposition \ref{prop:toriprop} is also true for centreless Lie tori
(in fact one can again take $\ell = 2$).  Hence, Theorem \ref{thm:realtorus}
is also true in that case.  That is, a centreless Lie torus $\cL$
over an algebraically closed field of characteristic 0
has a multiloop realization based on a finite dimensional  simple
algebra if and only if $\cL$ is fgc.  However, the analysis
similar to Remark \ref{rem:realtorus} is more subtle and requires
a more detailed study of Lie tori.  We therefore
postpone further discussion of this topic to a sequel to this paper.
\end{remark}

\subsection{A tensor product decomposition}
\label{subsec:tensorproduct}

\newcommand\cQ{\mathcal Q}

We now describe a tensor product decomposition for fgc associative tori.
If $m$ is a positive integer and $e$ is an integer that  is relatively prime to $m$,
we use the notation $\cQ(m,e)$ for the
quantum torus determined by the $2\times 2$-matrix
$\left[\begin{smallmatrix} 1&\zeta_m^e\\ \zeta_m^{-e} & 1 \end{smallmatrix}\right]$.
$\cQ(m,e)$ with its natural $\bbZ^2$-grading (see Example \ref{ex:quantum}) is an associative $\bbZ^2$-torus.
Also $k[z_1^{\pm 1},\dots, z_s^{\pm 1}]$ denotes the algebra of Laurent polynomials with its
natural $\bbZ^s$-grading.

\begin{proposition}
\label{prop:tensordecomp}
Suppose that $k$ is  algebraically closed of characteristic $0$ and $\Lm$ is a free abelian
group of finite rank $n\ge 1$.  Let $\cB$ be an fgc associative $\Lm$-torus.  Then
\begin{equation}
\label{eq:tensordecomp}
\cB \ig \cQ(m_1,e_1) \otk \dots \otk \cQ(m_r,e_r) \otk k[z_1^{\pm 1},\dots, z_s^{\pm 1}],
\end{equation}
where $r\ge 0$, $s\ge 0$, $n = 2r+s$,
$m_1,\dots,m_r$ are integers $\ge 2$  such that $m_1\mid m_2 \mid \dots \mid m_r$,
and $e_1,\dots,e_r$ are integers so that $\gcd(e_i,m_i) = 1$  for all $i$.
(On the right hand side of \eqref{eq:tensordecomp},
the tensor product has the natural
$\smash{\underbrace{\Ztwo\oplus \dots \oplus \Ztwo}_r} \oplus \bbZ^s$-grading determined
by the gradings on the components.)
\end{proposition}

This proposition has been proved by  K.-H.~Neeb \cite[Theorem III.3]{Neeb} using a
normal form for skew-symmetric integral matrices and some additional
arguments linking that normal form to the proposition  (see Remark \ref{rem:Neeb} below).
In order to illustrate how our results can be used to deduce
results about infinite dimensional graded algebras from results
about finite dimensional graded algebras, we outline an alternate proof of
the proposition using the Correspondence Theorem. (We use the Correspondence
Theorem rather than the Realization Theorem, since the coordinate
free point of view is more  convenient here.)

\begin{proof}[Outline of a proof] Suppose that $k$, $\Lm$ and $\cB$ are as in
the proposition.  Since $\supp(\cB) = \Lm$, it follows from
Proposition \ref{prop:bornmodfree} that $\Lm/\Gm(\cB)$ is finite. Let
$\bLm = \Lm/\Gm(\cB)$ and let $\pi : \Lm \to \bLm$ be the natural
projection. Now by the Correspondence Theorem, $\cB\isomLm L_\pi(\cA)$ for some finite
dimensional central-simple algebra $\cA$ that is $\bLm$-graded. In
fact  $\cA$ is an associative $\bLm$-torus.  But, since $k$ is
algebraically closed, we can identify $\cA$ with the algebra
$M_\ell(k)$ of $\ell\times\ell$-matrices over $k$ for some
$\ell\ge 1$. Furthermore, gradings on $M_\ell(k)$ by abelian
groups have been classified in \cite{BSeZ}.  Moreover, those for
which $M_\ell(k)$ is a $\bLm$-torus have a particularly simple
description \cite[Theorem 5]{BSeZ}, which gives a tensor product
decomposition for $\cA$ analogous to \eqref{eq:tensordecomp}
(without the divisibility condition  $m_1\mid m_2 \mid \dots \mid
m_r$). One can easily adjust this decomposition of $\cA$ to get
the divisibility condition. Then, in the final step of the proof
one shows (arguing as in the fundamental theorem
for finitely generated abelian groups)
that the decomposition of $\cA$
determines a  decomposition of the loop algebra $\loopA$.
\end{proof}

\begin{remark}
\label{rem:Neeb}
In \cite{Neeb},  Neeb proves more than is stated above.  He shows that
in \eqref{eq:tensordecomp} one can choose all $e_i = 1$ if $s>0$
and one can choose $e_i = 1$ for $2\le i \le r$ if $s = 0$.
This sharper result can be deduced from the one stated above
by means of a change of variables (using the argument in the proof of \cite[Theorem III.1]{Neeb}).
We note also that Neeb's  result in \cite{Neeb} is formulated in such a way that it holds
over any base  field $k$.
\end{remark}

\begin{remark}
Alternative  tori have been classified in general
in \cite{BGKN}. Although we have not worked out the details,
it should be straightforward to alternately apply
our approach above
to obtain a  classification  of all fgc
alternative tori.  This would use the description by
A.~Elduque of all gradings of the octonion algebra over $k$~\cite{E}.

Jordan tori  have been classified in general
by Yoshii in \cite{Y1}.   Yoshii used deep results
on prime Jordan algebras due to Zelmanov, and so a more elementary
approach would also be interesting.
Our approach to classifying  fgc Jordan tori would  require a description
of the gradings on finite dimensional simple Jordan algebras
that satisfy conditions (a) and (b) in Remark
\ref{rem:realtorus}.   Such a description
does not seem yet to have been completed, although
a lot is known about finite dimensional gradings
(see for example \cite{BShZ} and the references  therein).
\end{remark}

\end{document}